\documentclass[a4paper, 11pt]{article}

\usepackage[utf8]{inputenc}
\usepackage[english]{babel}

\usepackage{amsmath,amssymb,amsfonts}
\usepackage{amsthm}
\usepackage{bm}
\usepackage{dsfont}

\usepackage[paperwidth=192mm,
		    		   paperheight=262mm,
		    		   vmargin={19mm,19mm},
		    		   hmargin={13.7mm,13.7mm},
		    		   headsep=12pt,
		    		   footskip=12pt]{geometry}

\usepackage{siunitx}
\usepackage{hyperref}

\usepackage{float}

\usepackage{caption}
\usepackage{subcaption}
\usepackage{booktabs}
\usepackage{tabularx}

\usepackage{pgfplots}
\pgfplotsset{compat=1.16}

\usepackage{tikz}
\usetikzlibrary{shapes,positioning,intersections,quotes}
\tikzset{
    cross/.pic = 
    {
        \draw[rotate = 45] (-#1,0) -- (#1,0);
        \draw[rotate = 45] (0,-#1) -- (0, #1);
    }
}

\usepackage[backend=bibtex,
					   style=numeric,
            		   sorting=nyt,
            		   defernumbers=true]{biblatex}
\renewbibmacro{in:}{}
\DeclareBibliographyCategory{cited}
\AtEveryCitekey{\addtocategory{cited}{\thefield{entrykey}}}
\DeclareMathOperator{\grad}{\nabla}
\DeclareMathOperator{\dive}{\nabla\cdot}

\newcommand{\vel}{\bm{u}}

\newcommand{\pad}[2]{\frac{\partial{#1}}{\partial{#2}}}

\newcommand{\rpth}[1]{\left(#1\right)}
\newcommand{\spth}[1]{\left[#1\right]}
\newcommand{\cpth}[1]{\left\{#1\right\}}

\newtheorem{theorem}{Theorem}[section]
\newtheorem{definition}[theorem]{Definition}

\nocite{*}

\begin{document}

\pagenumbering{arabic}

\title{Second-order optimally stable IMEX (pseudo-)staggered Galerkin discretization with application to depth-integrated lava flow simulations}

\date{}

\author{Federico Gatti$^{(1)}$, Giuseppe Orlando$^{(2)}$}

\maketitle

\begin{center}
{\small
$^{(1)}$  
Seminar for Applied Mathematics (SAM),\\ Department of Mathematics, ETH Z\"{u}rich,\\ CH-8092 Z\"{u}rich, Switzerland \\
{\tt federico.gatti@math.ethz.ch}
\vskip 0.2cm
$^{(2)}$  
CMAP, CNRS, \'{E}cole polytechnique, Institut Polytechnique de Paris \\ 
Route de Saclay, 91120 Palaiseau, France \\
{\tt giuseppe.orlando@polytechnique.edu}
}
\end{center}

\noindent


%%%%%%%%%%%%%%%% Abstract %%%%%%%%%%%%%%%%%%
\begin{abstract} \noindent
We present second-order optimally stable Implicit-Explicit (IMEX) Runge--Kutta (RK) schemes with application to a modified set of shallow water equations that can be used to model the dynamics of lava flows. The schemes are optimally stable in the sense that they satisfy, at the space-time discretization level, a condition analogous to the \texttt{L}-stability of Runge--Kutta methods for ordinary differential equations. A novel (pseudo-)staggered Galerkin scheme is introduced, which can be interpreted as an extension of the classical two-step Taylor--Galerkin (TG2) scheme. The method is derived by combining a von Neumann stability analysis with a Lax--Wendroff procedure. For the discretization of the non-conservative terms that characterize the lava flow model, we employ the Path-Conservative (PC) method. The proposed scheme is evaluated on a number of relevant test cases, demonstrating accuracy, robustness, and well-balancing properties for the lava flow model.
\end{abstract}
\setlength{\leftskip}{0pt}
\setlength{\rightskip}{0pt}	

\pagebreak

%%%%%%%%%% Introduction %%%%%%%%%%%%%%%%%%%%%
\section{Introduction}
\label{sec:intro} \indent

Volcanic eruptions represent a natural phenomenon that threatens millions of people around the world. Eruptions may show effusive activity with the propagation of a lava flow from the vent~\cite{biagioli:2023, gatti:2025}. The risks and damages associated with the propagation of lava flows require a mathematical description of this phenomenon and the ability to perform numerical simulations, to improve forecast skills and decision-making in volcanic risk management~\cite{costa:2005, gatti:2025}. A widely used mathematical description of lava flows is based on a depth-integrated model, originally proposed in~\cite{costa:2005}, and then extended in~\cite{biagioli:2021, biagioli:2023, conroy:2021, gatti:2025}. The model is a modification of the classical shallow water equations (SWE)~\cite{decoene:2009}, also known as the Saint-Venant (or de Saint-Venant) equations. SWE are a simplified model for hydrostatic flows, in which the depth of the fluid is much smaller than its width, and are widely used in river hydraulics, coastal engineering, flood problems, and landslide runout forecasting; see, among many others,~\cite{arpaia:2026, bresch:2009, ciallella:2022, gatti:2023a, gatti:2023b, gatti:2024b, giraldo:2002, guermond:2025, kinnmark:2012, kurganov:2007, soares:2008, williamson:1992} and the references therein. The lava flow model includes the depth-integration of the energy equation from the complete Navier--Stokes system~\cite{costa:2005}. Temperature plays a significant role in the dynamics of lava flows and cannot be neglected. More specifically, the presence of vents is modeled as a point source in the depth-integrated model~\cite{gatti:2025}.

Considering lava flow simulations as a numerically challenging testbed, the goal of this paper is to present a novel (pseudo-)staggered Galerkin method that can be interpreted as an extension of the classical two-step Taylor--Galerkin (TG2) scheme~\cite{donea:1984}. The TG2 scheme has been used for the solution of a single-phase single-layer depth-integrated landslide model~\cite{gatti:2024a}, of a two-phase two-layer landslide model~\cite{gatti:2024b}, as well as the lava flow model~\cite{gatti:2025}.

The lava flow model contains non-conservative terms. We employ the Path-Conservative (PC) method~\cite{busto:2021, castro:2009, pares:2006} for their discretization as done, e.g., in~\cite{gatti:2024b, gatti:2024a, gatti:2025}. The time discretization is based on an Implicit-Explicit (IMEX) Runge--Kutta (RK) method. Since the seminal papers~\cite{casulli:1990, casulli:1992}, semi-implicit time discretization methods are extremely popular for SWE and have been extended to high-order finite elements, e.g., in~\cite{boscheri:2023, busto:2022, dumbser:2013, tavelli:2014}. IMEX methods have been used for the lava flow model in~\cite{biagioli:2023, gatti:2025}. The presence of source terms related to friction in the lava flow model implies a potential source of stiffness. Hence, a popular strategy consists in implicitly treating these source terms, while the other terms are treated in an explicit fashion~\cite{gatti:2025}.

IMEX methods are typically derived in the Method of Lines (MOL) framework, and this is also the philosophy adopted, e.g., in~\cite{boscheri:2023} for SWE. However, recent works on RK-DG methods have shown that RK time integration may not preserve important stability properties of the semi-discrete scheme, such as $L^{2}$ or entropy stability~\cite{ranocha:2020, sun:2019, xu:2019}. For this reason, we derive our IMEX scheme slightly beyond the traditional MOL approach. More specifically, we analyze the amplification factor obtained from a von Neumann stability analysis~\cite{deriaz:2020}, in which both space and time discretizations are taken into account. For this purpose, we introduce the novel concept of space-time \texttt{L}-stability. It has been shown~\cite{gatti:2025} that oscillations in mass flux can arise if the amplification factor does not exhibit the desired decay properties. In that earlier work, the coefficients of the TG2 scheme have been adopted as an explicit companion of an IMEX scheme (see also~\cite{liu:2006}). Here, instead, we derive new IMEX coefficients designed to enforce improved decay properties (see Section~\ref{ssec:staggered_von_Neumann_et_al}). As in~\cite{gatti:2025}, intermediate solutions are defined in the discrete space $\bm{Q}_{0}$ on the primal grid, while the final stage/solution belongs to $\tilde{\bm{Q}}_{1}$ (see Section~\ref{sec:space_disc}), where the mass grouping relocates the degrees of freedom to the dual grid. Hence, the method can be interpreted both as an extension of TG2 and as a (pseudo-)staggered Galerkin scheme, since intermediate and final stages/solutions effectively live on different grids. This choice does not reduce the overall order of the scheme as one can notice from the discretization error arising from the modified equation analysis (see Section~\ref{ssec:staggered_von_Neumann_et_al}). Indeed, intermediate stages do not necessarily correspond to fully consistent spatial approximations of the final solution, but serve as internal states in the IMEX time integration process. A Lax--Wendroff procedure is finally used to complete the derivation of a scheme with enhanced stability properties.

The paper is structured as follows. In Section~\ref{sec:model}, we discuss the set of equations of the lava flow model. In Section~\ref{sec:time_disc}, we present the IMEX time discretization method and describe the characteristics of the proposed IMEX scheme. In Section~\ref{sec:space_disc}, we recall the space discretization strategy and some concepts related to the Path-Conservative method. We present a detailed derivation of the numerical method through a von Neumann stability analysis and a Lax--Wendroff procedure, and we also show that the method is of the second order in spite of the use of the (pseudo-)staggered approach. Some numerical results to assess the robustness, accuracy and well-balancing properties for the lava flow model of this novel scheme are presented in Section~\ref{sec:num}. Finally, some conclusions and perspectives for future work are discussed in Section~\ref{sec:conclu}.

%%%%%%%%%%% Mathematical model %%%%%%%%%%%%
\section{The mathematical model}
\label{sec:model}

The mathematical model consists of a set of equations obtained from the depth-integration of the Navier--Stokes equations under the hypothesis of hydrostatic pressure distribution along the vertical axis~\cite{biagioli:2021, costa:2005}. Let $\Omega = \rpth{0, L_{x}} \times \rpth{0, L_{y}}$ be a Cartesian domain and denote by $\bm{x}$ the spatial coordinates and by $t$ the temporal coordinate. We define the space-time dependent domain $\Omega_{w} \subset \Omega$ as the subdomain where the height of the lava $h$ is greater than zero. The model equations, endowed with suitable initial and boundary conditions, read as follows
\begin{align}\label{eq:lava_model_Z}
	\pad{h}{t} + \dive\rpth{h\vel} &= Q\delta\rpth{\bm{x} - \bm{x}_{v}}, \nonumber \\
	\pad{(h\vel)}{t} + \dive\rpth{h\vel \otimes \vel} + \grad\rpth{\frac{1}{2}gh^{2}} + gh\grad Z &= -\lambda\rpth{h, T}\vel, \\
	\pad{(hT)}{t} + \dive\rpth{hT\vel} &= Q T_{e}\delta\rpth{\bm{x} - \bm{x}_{v}}, \nonumber
\end{align}
for $\bm{x} \in \Omega_{w}, t \in \left(0, T_{f}\right]$, with $T_{f}$ being the final time. Here, $\otimes$ denotes the tensor product, $g$ is the Earth acceleration of gravity, $\vel = \rpth{u_{x}, u_{y}}^{\top}$ is the vector of vertically averaged velocities, $T$ is the temperature integrated with depth and $Z = Z(\bm{x})$ is the bottom topography. For simplicity, we assume that $Z$ does not change over time. Moreover, $\delta(\cdot)$ denotes the Dirac delta distribution, $\bm{x}_{v}$ is the spatial location of the vent, $Q = Q(t)$ is the lava discharge into the domain, which may vary in time, and $T_{e}$ is the vent effusion temperature. Finally, $-\lambda\vel$ represents a friction term, for which the positive coefficient $\lambda$ has the following expression
\begin{equation}\label{eq:lambda_vent}
	\lambda = \frac{3\nu_{r}}{h}e^{-b\rpth{T - T_{r}}},
\end{equation}
where $\nu_{r}$ is a reference kinematic viscosity, $T_{r}$ is a reference temperature, and $b$ is a suitable coefficient. All variables and coefficients are in SI units. Following~\cite{costa:2005}, we choose the reference temperature $T_{r}$ equal to the vent effusion temperature $T_{e}$.

By rearranging the momentum equation in system~\eqref{eq:lava_model_Z} and explicitly differentiating the pressure term, the governing equations of the lava flow model can be reformulated as follows~\cite{arpaia:2026, tavelli:2014}
\begin{align}\label{eq:lava_model_non_conservative}
	\pad{h}{t} + \dive\rpth{h\vel} &= Q\delta\rpth{\bm{x} - \bm{x}_{v}}, \nonumber \\
	\pad{(h\vel)}{t} + \dive\rpth{h\vel \otimes \vel} + gh\grad\zeta &= -\lambda\rpth{h, T}\vel, \\
	\pad{(hT)}{t} + \dive\rpth{hT\vel} &= Q T_{e}\delta\rpth{\bm{x} - \bm{x}_{v}}, \nonumber
\end{align}
where $\zeta = h + Z$ is the total height of the free surface. One can easily notice that, without vents, systems~\eqref{eq:lava_model_Z} and~\eqref{eq:lava_model_non_conservative} admit the following steady-state solution
\begin{equation}
	\zeta = h + Z = \text{constant}, \qquad \vel = \bm{0}, \qquad T = \text{constant},
\end{equation}
of the so-called lake-at-rest condition. In Section~\ref{sec:num}, we will analyze the ability of the presented numerical method to preserve this condition, the so-called well-balancing property~\cite{castro:2013}, also known as {\cal C}-property~\cite{vazquez:1999}.

%%%%%%%%%%% Time discretization %%%%%%%%%%%
\section{The time discretization method}
\label{sec:time_disc}

In this section, we briefly outline the time discretization strategy. An Implicit-Explicit (IMEX) Runge--Kutta (RK) method \cite{boscarino:2024, kennedy:2003} is employed. These methods are rather popular for time dependent problems that can be formulated as $\mathbf{y}' = \mathbf{f}_{\mathrm{S}}(\mathbf{y},t) + \mathbf{f}_{\mathrm{NS}}(\mathbf{y},t)$, where the $\mathrm{S}$ and $\mathrm{NS}$ subscripts denote the stiff and non-stiff components of the system, to which the implicit and explicit companion methods are applied, respectively. IMEX methods are represented compactly by the following two Butcher tableaux~\cite{butcher:2008}
\begin{center}
	\begin{tabular}{c|c}
		$\mathbf{c}$ & $\mathbf{A}$ \\
		\hline
		\rule{0pt}{2.5ex} & $\mathbf{b}^{\top}$
	\end{tabular}
	\qquad
	\begin{tabular}{c|c}
		$\tilde{\mathbf{c}}$ & $\tilde{\mathbf{A}}$ \\
		\hline
		\rule{0pt}{2.5ex} & $\tilde{\mathbf{b}}^{\top}$
	\end{tabular}
\end{center}
with $\mathbf{A} = \cpth{a_{lm}}, \mathbf{b} = \cpth{b_{l}}, \mathbf{c} = \cpth{c_{l}}, \tilde{\mathbf{A}} = \cpth{\tilde{a}_{lm}}, \tilde{\mathbf{b}} = \cpth{\tilde{b}_{l}}$, and $\tilde{\mathbf{c}} = \cpth{\tilde{c}_{l}}$. Coefficients $a_{lm}, \tilde{a}_{lm}, c_{l}, \tilde{c}_{l}, b_{l}$, and $\tilde{b}_{l}$ are determined so that the method is consistent of a given order. In particular, the following relation has to be satisfied for the method to be consistent~\cite{hairer:1996, kennedy:2003}
\begin{equation}\label{eq:order1_condition}
	\sum_{l=1}^{s}b_{l} = \sum_{l=1}^{s}\tilde{b}_{l} = 1.
\end{equation}
We also assume that
\begin{equation}\label{eq:compatibility_condition}
	\sum_{m=1}^{s}a_{lm} = c_{l}, \qquad \sum_{m=1}^{s}\tilde{a}_{lm} = \tilde{c}_{l},
\end{equation}
which is a standard assumption for Runge--Kutta schemes. This assumption simplifies the order conditions and guarantees that each stage achieves at least first-order accuracy~\cite{ketcheson:2009, wanner:1996}. Moreover, for non-autonomous problems, condition~\eqref{eq:compatibility_condition} ensures that the method produces the same result whether the problem is formulated in autonomous or non-autonomous form~\cite{ketcheson:2023}.

We adopt a Diagonally Implicit Runge--Kutta (DIRK) scheme for the implicit part, which means that $\tilde{a}_{lm} = 0$ for $l < m$. If $\mathbf{v}^{n} \approx \mathbf{y}(t^{n})$, the generic $l$-th stage of the IMEX method can be written as follows
\begin{align}
	\mathbf{v}^{(n,l)} &= \mathbf{v}^{n} + 
	\Delta t \sum_{m=1}^{l-1} a_{lm}\mathbf{f}_{NS}\rpth{\mathbf{v}^{(n,m)}, t + c_{m} \Delta t} \nonumber \\
	&+ \Delta t \sum_{m=1}^{l-1} \tilde{a}_{lm}\mathbf{f}_{S}\rpth{\mathbf{v}^{(n,m)}, t + \tilde{c}_{m} \Delta t} + \Delta t\tilde{a}_{ll}\mathbf{f}_{S}\rpth{\mathbf{v}^{(n,l)}, t + \tilde{c}_{l}\Delta t},
\end{align}
where $l=1, \dots, s$, with $s$ denoting the number of stages, and $\Delta t$ denotes the time step. After the computation of all intermediate stages, $\mathbf{v}^{n+1}$ is updated as follows
\begin{equation}\label{eq:update_IMEX}
	\mathbf{v}^{n+1} = \mathbf{v}^{n} + \Delta t \sum_{l=1}^{s} b_{l} \mathbf{f}_{NS}\rpth{\mathbf{v}^{(n,l)}, t + c_{l}\Delta t} + \Delta t \sum_{l=1}^{s} \tilde{b}_{l} \mathbf{f}_{S}\rpth{\mathbf{v}^{(n,l)}, t + \tilde{c}_{l}\Delta t}.
\end{equation}

Hence, the generic $l$-th IMEX stage applied to system~\eqref{eq:lava_model_Z}, in which the stiff friction source term is treated implicitly, while the remaining terms are treated explicitly, reads
\begin{align}\label{eq:imex_lava_model_Z}
	h^{(n,l)} &= h^{n} - \Delta t\sum_{m=1}^{l-1}a_{lm}\dive\rpth{\rpth{h\vel}^{(n,m)}} + \Delta t\sum_{m=1}^{l-1}{a}_{lm}Q\rpth{t^{n} + {c}_{m}\Delta t}\delta\rpth{\bm{x} - \bm{x}_{v}}, \nonumber \\
	\rpth{h\vel}^{(n,l)} &= \rpth{h\vel}^{n} - \Delta t\sum_{m=1}^{l-1}a_{lm}\dive\rpth{\rpth{h\vel}^{(n,m)} \otimes \vel^{(n,m)}} \nonumber \\
	&- \Delta t\sum_{m=1}^{l-1}a_{lm}\grad\rpth{\frac{1}{2}g\rpth{h^{(n,m)}}^{2}} \\
	&- \Delta t\sum_{m=1}^{l-1}\tilde{a}_{lm}\lambda\rpth{h^{(n,m)}, T^{(n,m)}}\vel^{(n,m)} - \Delta t \tilde{a}_{ll}\lambda\rpth{h^{(n,l)}, T^{(n,l)}}\vel^{(n,l)}, \nonumber \\
	\rpth{hT}^{(n,l)} &= \rpth{hT}^{n} - \Delta t\sum_{m=1}^{l-1}a_{lm} \dive\rpth{\rpth{h\vel}^{(n,m)}T^{(n,m)}} \nonumber \\
	&+ \Delta t\sum_{m=1}^{l-1}{a}_{lm}Q\rpth{t^{n} + {c}_{m}\Delta t}\delta\rpth{\bm{x} - \bm{x}_{v}}. \nonumber
\end{align}
It is worth noting that source terms independent of the state variables, i.e. those associated with vent modeling, are integrated using the explicit Butcher tableau.
The solution update then reads
\begin{align}\label{eq:update_lava_model}
	h^{n+1} &= h^{n} - \Delta t\sum_{l=1}^{s}b_{l}\dive\rpth{\rpth{h\vel}^{(n,l)}} + \Delta t\sum_{l=1}^{s}{b}_{l}Q\rpth{t^{n} + {c}_{l}\Delta t}\delta\rpth{\bm{x} - \bm{x}_{v}}, \nonumber \\
	\rpth{h\vel}^{n+1} &= \rpth{h\vel}^{n} - \Delta t\sum_{l=1}^{s}b_{l}\dive\rpth{\rpth{h\vel}^{(n,l)} \otimes \vel^{(n,l)}} \nonumber \\
	&- \Delta t\sum_{l=1}^{s}b_{l}\grad\rpth{\frac{1}{2}g\rpth{h^{(n,l)}}^{2}} 
	- \Delta t\sum_{l=1}^{s}\tilde{b}_{l}\lambda\rpth{h^{(n,l)}, T^{(n,l)}}\vel^{(n,l)}, \\
	\rpth{hT}^{n+1} &= \rpth{hT}^{n} - \Delta t\sum_{l=1}^{s}b_{l}\dive\rpth{\rpth{h\vel}^{(n,l)}T^{(n,l)}} \nonumber \\
	&+ \Delta t\sum_{l=1}^{s}{b}_{l}Q\rpth{t^{n} + {c}_{l}\Delta t}\delta\rpth{\bm{x} - \bm{x}_{v}}. \nonumber
\end{align}

In this work, we focus on three-stage ($s = 3$) second-order IMEX methods. Higher-order IMEX methods can be found, e.g., in~\cite{orlando:2025a, orlando:2025b}. We now provide an overview of the characteristic of the proposed IMEX scheme. The main novelty of this work, as will become clear later, is that some of the coefficients are determined by taking the spatial discretization into account, rather than being treated as if one were solving a semi-discrete system of ordinary differential equations (ODEs), as in the classical MOL framework. This approach also incorporates a novel notion of stability, namely the space-time \texttt{L}-stability, which will be introduced in Section~\ref{ssec:staggered_von_Neumann_et_al}.

\subsection{Characteristics of the proposed IMEX scheme}
\label{ssec:novel_IMEX}

In this section, the main properties of the proposed IMEX scheme are discussed. The order and coupling conditions for second-order IMEX methods read~\cite{pareschi:2005, rokhzadi:2018}
%.
\begin{equation}\label{eq:coupling_conditions}
	\sum_{i=1}^{s}b_{i}c_{i} = \frac{1}{2}, \qquad \sum_{i=1}^{s} \tilde{b}_{i}\tilde{c}_{i} = \frac{1}{2}, \qquad \sum_{i=1}^{s}\tilde{b}_{i}c_{i} = \frac{1}{2}, \qquad \sum_{i=1}^{s}b_{i}\tilde{c}_{i} = \frac{1}{2}.
\end{equation}
We consider for the implicit companion method the family of the so-called ESDIRK (Explicitly Singly Diagonal Implicit Runge--Kutta) integration methods~\cite{kennedy:2016}, for which $\tilde{a}_{11} = 0$ and $\tilde{a}_{ll} = \gamma > 0, 2 \le l \le s$. Hence, the first stage of the IMEX scheme~\eqref{eq:imex_lava_model_Z} is purely formal and reads $\mathbf{q}^{(n,1)} = \mathbf{q}^{n}$. Following, e.g.,~\cite{boscarino:2007, giraldo:2013, rokhzadi:2018}, we employ a stiffly accurate (SA) scheme for the implicit method, which is responsible of the integration of stiff source terms. Hence, $\tilde{a}_{sk} = b_{k}, k = 1,\dots,s$. This property also leads to computational efficiency for the implicit method, since the output of the last stage is equal to the updated solution. For what concerns the explicit companion method, we assume $\mathbf{b} = \tilde{\mathbf{b}}$, which simplifies the order conditions and, more importantly, it is a necessary condition for the preservation of linear invariants~\cite{giraldo:2013}. Hence, after application of~\eqref{eq:order1_condition},~\eqref{eq:compatibility_condition}, and~\eqref{eq:coupling_conditions}, the Butcher tableaux read as follows

\begin{table}[H]
	\begin{minipage}{0.45\textwidth}
		\begin{center}
			\begin{tabular}{c|ccc}
				$0$ & $0$ & $0$ & $0$ \\
				$a_{21}$ & $a_{21}$ & $0$ & $0$ \\
				$a_{31} + a_{32}$ & $a_{31}$ & $a_{32}$ & $0$ \\
				\hline
				& $\tilde{a}_{31}$ & $\tilde{a}_{32}$ & $\gamma$
			\end{tabular}
		\end{center}
	\end{minipage}\hspace{1cm}
	\begin{minipage}{0.45\textwidth}
		\begin{center}
			\begin{tabular}{c|ccc}
				$0$ & $0$ & $0$ & $0$ \\
				$\tilde{a}_{21}+\gamma$ & $\tilde{a}_{21}$ & $\gamma$ & $0$ \\
				1 & $\tilde{a}_{31}$ & $\tilde{a}_{32}$ & $\gamma$ \\
				\hline
				& $\tilde{a}_{31}$ & $\tilde{a}_{32}$ & $\gamma$
			\end{tabular}
		\end{center}
	\end{minipage}
	\caption{Butcher tableaux of the the IMEX-RK2 method. Left: explicit method. Right: implicit method.}
	\label{tab:imex_rk2_butch}
\end{table}
\noindent
where
\begin{subequations}\label{eq:number_1_cond}
	\begin{align}
		\tilde{a}_{31} &= 1 - \gamma - \frac{1 - 2\gamma}{2\rpth{\tilde{a}_{21} + \gamma}} \qquad \text{(consistency condition \eqref{eq:order1_condition})}, \\
		\tilde{a}_{32} &= \frac{1 - 2\gamma}{2\rpth{\tilde{a}_{21} + \gamma}} \qquad \text{(order/coupling condition $\sum_{i}\tilde{b}_{i}\tilde{c}_{i} = \frac{1}{2}$)}, \label{eq:a32_tilde}
	\end{align}
\end{subequations}
and $a_{21}$ is such that
\begin{equation}\label{eq:number_2_cond}
	\frac{1 - 2\gamma}{2\gamma\rpth{\tilde{a}_{21} + \gamma}}a_{21} + \gamma\rpth{a_{31} + a_{32}} = \frac{1}{2} \quad \text{(order/coupling condition $\sum_{i}\tilde{b}_{i}c_{i} = \frac{1}{2}$)}.
\end{equation}
Indeed, since $\mathbf{b} = \tilde{\mathbf{b}}$, only two order and coupling conditions~\eqref{eq:coupling_conditions} apply, and we are therefore left with four free parameters.

A minimal requirement for the implicit method is to be \texttt{A}-stable. A Runge--Kutta method is called \texttt{A}-stable if $\left|R(z)\right| \le 1$ for $Re(z) \le 0$~\cite{kennedy:2016}, where $R(z)$ denotes the stability function. or ESDIRK methods the stability function can be written as~\cite{boscarino:2007, hairer:1996, kennedy:2016}
\begin{equation}\label{eq:stab_function_RK}
	R(z) = \frac{\mathrm{\det}\rpth{\mathbf{I} - z\tilde{\mathbf{A}} + z\mathds{1}\tilde{\mathbf{b}}^{\top}}}{\mathrm{\det}\rpth{\mathbf{I} - z\tilde{\mathbf{A}}}} = \frac{P(z)}{Q(z)} = \frac{P(z)}{\rpth{1 - \gamma z}^{s-1}},
\end{equation}
where $\mathbf{I} \in \mathbb{R}^{s \times s}$ is the identity matrix and $\mathds{1} = \rpth{1, 1, \dots, 1}^{\top}$. $P(z)$ is a polynomial of the highest degree $s-1$. Since $\gamma > 0$, the stability function $R(z)$ is analytical in the complex left half-plane and therefore does not contain any pole. Hence, following the result in~\cite{alt:1972}, if the method is \texttt{I}-stable, i.e., $\left|R(z)\right| \le 1$ for $Re(z) = 0$, then it is also \texttt{A}-stable. In order to test the \texttt{I}-stability, we consider the \texttt{E}-polynomial introduced in~\cite{norsett:1975}
\begin{equation}
	E(y) = \left|Q(iy)\right|^{2} - \left|P(iy)\right|^{2} \qquad y \in \mathbb{R}.
\end{equation}
\texttt{I}-stability requires that $E(y) \ge 0 \quad \forall y \in \mathbb{R}.$ By following the procedure described in \cite[pp. 40-41]{hairer:1996}, one can verify that the stability function~\eqref{eq:stab_function_RK} of the implicit companion method in Table~\ref{tab:imex_rk2_butch} with \eqref{eq:number_1_cond} reduces to
\begin{equation}
	R(z) = \frac{\rpth{1 - \gamma z}^{2} - 2\gamma z^{2} + \frac{z^{2}}{2} + z}{\rpth{1 - \gamma z}^{2}},
\end{equation}
so that, after some algebraic manipulations,
\begin{equation}
	E(y) = \rpth{4\gamma^{3} - 5\gamma^{2} + 2\gamma - \frac{1}{4}}y^{4}.
\end{equation}
The implicit method is therefore \texttt{I}-stable (hence \texttt{A}-stable) if $\gamma \ge \frac{1}{4}$.

System~\eqref{eq:lava_model_Z}-\eqref{eq:lava_model_non_conservative} tends to an index-1 differential-algebraic system as $\lambda \to \infty$~\cite{boscarino:2007, boscarino:2009}. In this case, supplementary conditions are required to ensure that the order convergence is preserved for the differential-algebraic system. More specifically, we require~\cite{boscarino:2009}
\begin{equation}
	\hat{\mathbf{b}}^{\top}\hat{\boldsymbol{\mathcal{A}}}^{-1}\hat{\mathbf{c}} = 1,
\end{equation}
where $\hat{\mathbf{b}}^{\top} = \spth{\tilde{b}_{l}}_{l>1}, \hat{\boldsymbol{\mathcal{A}}} = \spth{\tilde{a}_{lm}}_{l,m>1},$ and $\hat{\mathbf{c}} = \spth{\tilde{c}_{l}}_{l>1}$. One can easily verify that this condition is fulfilled for any value of $\gamma$ for the implicit method in Table~\ref{tab:imex_rk2_butch}.

Another important property of implicit time discretization methods is the so-called \texttt{L}-stability. A Runge--Kutta method is said to be \texttt{L}-stable~\cite{wanner:1996} if it is \texttt{A}-stable and $R(z) \to 0$ as $z \to \infty$. Since the seminal paper~\cite{pareschi:2005}, applying \texttt{L}-stable methods in a stiff nonlinear regime is widely regarded as the appropriate choice, because it provides the resulting method with an intrinsically robust response to high frequency disturbances that may arise. The concept of \texttt{L}-stability is introduced for the analysis of ODE methods applied to linear constant-coefficient model problems. Following the result in~\cite{wanner:1996}, a \texttt{L}-stable method is typically obtained from the combination of an \texttt{A}-stable method with a SA scheme. However, this combination does not necessarily results into a \texttt{L}-stable method for Runge--Kutta methods for which the matrix $\tilde{\mathbf{A}}$ is not invertible~\cite{boscarino:2009}, as it occurs for ESDIRK methods. Moreover, as discussed in~\cite{gatti:2025}, oscillations in mass flux can arise if the amplification factor obtained from a von Neumann stability analysis does not satisfy decaying properties that are, to some extent, analogous to the \texttt{L}-stability property. Hence, we employ the amplification factor of a von Neumann stability analysis to determine the remaining coefficients of the implicit companion method, i.e. $\tilde{a}_{21}$ and $\gamma$, and we will further discuss in the following the link with the classical \texttt{L}-stability property. The von Neumann stability analysis provides also one parameter of the explicit companion method because of the (pseudo-)staggered Galerkin approach (see Section~\ref{ssec:staggered_von_Neumann_et_al}). The remaining free parameter is determined employing a Lax--Wendroff procedure so as to further enhance the stability properties of the resulting numerical method. We refer therefore to Section~\ref{ssec:staggered_von_Neumann_et_al} for the completion of the IMEX scheme. Finally, it is worth to notice that we are allowing $\mathbf{c} \neq \mathbf{\tilde{c}}$. The condition $\mathbf{c} = \mathbf{\tilde{c}}$ is usually employed to simplify the order conditions and the design of IMEX methods~\cite{kennedy:2003, kennedy:2016}, especially for high-order methods~\cite{kennedy:2003, kennedy:2019}, but it is not necessary. We will further discuss this point in Section~\ref{ssec:staggered_von_Neumann_et_al}.

%%%%%%% Space discretization %%%%%%%%%%%%%%
\section{The space discretization method}
\label{sec:space_disc}

In this section, we outline the space discretization method. For this purpose, we rewrite system~\eqref{eq:lava_model_Z} in the compact form
\begin{equation}\label{eq:lava_model_Z_compact}
	\pad{\mathbf{q}}{t} + \dive\mathbf{F}\rpth{\mathbf{q}} + \mathbf{B}\rpth{\mathbf{q}}\grad\mathbf{q} = \mathbf{S}_{\mathrm{I}}\rpth{\mathbf{q}} + \mathbf{S}_{\mathrm{E}},
\end{equation}
\begin{subequations}
	\begin{align}
		&\mathbf{q} =
		\begin{bmatrix}
			h \\
			h\vel \\
			hT \\
			Z
		\end{bmatrix},
		\qquad
		\mathbf{F}\rpth{\mathbf{q}} =
		\begin{bmatrix}
			h\vel^{\top} \\
			h\vel \otimes \vel + \frac{1}{2}gh^{2}\mathbf{I} \\[1.5mm]
			hT\vel^{\top} \\[1.5mm]
			\mathbf{0}^{\top}
		\end{bmatrix} \\
		&\mathbf{S}_{\mathrm{I}}\rpth{\mathbf{q}} =
		\begin{bmatrix}
			0 \\
			-\lambda\rpth{h,T}\vel \\
			0 \\
			0
		\end{bmatrix},
		\qquad
		\mathbf{S}_{\mathrm{E}} =
		\begin{bmatrix}
			Q\delta\rpth{\bm{x} - \bm{x}_{v}} \\
			\bm{0} \\
			QT_{e}\delta\rpth{\bm{x} - \bm{x}_{v}} \\
			0
		\end{bmatrix},
	\end{align}
\end{subequations}
with $\mathbf{I}$ denoting the identity tensor. Moreover, $\mathbf{B}\rpth{\mathbf{q}}$ entails a block-matrix syntax so as to obtain a compact notation. More specifically, 
$$\mathbf{B}\rpth{\mathbf{q}} = \rpth{\mathbf{0},\mathbf{0},\mathbf{0},\hat{\mathbf{B}}\rpth{\mathbf{q}}},$$
where
\begin{equation}\label{eq:non_cons_B}
	\hat{\mathbf{B}}\rpth{\mathbf{q}} =
	\begin{bmatrix}
		\mathbf{0}^{\top} \\
		gh\mathbf{I} \\
		\mathbf{0}^{\top} \\
		\mathbf{0}^{\top}
	\end{bmatrix}.
\end{equation}
We consider a decomposition of the domain $\Omega$ into a family of quadrilaterals $\mathcal{T}_{\mathcal{H}}$ and denote each element by $K$. We introduce the following finite element spaces
\begin{align*}
	Q_{0} &= \cpth{v \in L^{2}(\Omega) : v\rvert_K \in \mathbb{Q}_{0} \quad \forall K \in \mathcal{T}_{\mathcal{H}}}, \qquad \bm{Q}_{0} = \spth{Q_{0}}^{d + 2}, \\
	\tilde{Q}_{1} &= \cpth{v \in C^{0}(\Omega) : v\rvert_K \in \mathbb{Q}_{1}, \quad \forall K \in \mathcal{T}_{\mathcal{H}}}, \qquad \tilde{\bm{Q}}_{1} = \spth{\tilde{Q}_{1}}^{d + 2},
\end{align*}
where $d$ denotes the number of spatial dimensions and $\mathbb{Q}_{r}, r = 0,1$ is the space of polynomials of degree $r$ in each coordinate direction. We seek the intermediate solution $\mathbf{q}^{(n,2)}$ in the discrete space $\bm{Q}_{0}$ and the solutions $\mathbf{q}^{(n,3)}$ and $\mathbf{q}^{n+1}$ in the discrete space $\tilde{\bm{Q}}_{1}$, the latter being obtained by performing a mass lumping of the mass matrices arising from the $L^{2}$ projection onto the bilinear space $\tilde{\bm{Q}}_{1}$. Notice that, similarly to the classical TG2 scheme, the scheme achieves second-order accuracy in space despite the use of piecewise constant elements for the intermediate stage~\cite{gatti:2023c, gatti:2023b}: the discretization error introduced by the piecewise constant ($\bm{Q}_{0}$) representation at this stage is compensated by the subsequent mass lumping and reconstruction steps, so that the overall scheme effectively recovers a second-order accurate flux evaluation. This mass lumping, applied in conjunction with the $\bm{Q}_{0}$--$\tilde{\bm{Q}}_{1}$ discretization, leads to a formulation that is structurally analogous to a staggered finite volume method. Once the time discretization has been introduced, the modified equation presented in Equation~\eqref{eq:modified_equation} will indeed confirm the second-order accuracy of the overall scheme.
The discrete weak form associated to the second stage of the IMEX scheme reads therefore as follows
\begin{align}
	\int_{\Omega}\mathbf{q}^{(n,2)} \cdot \boldsymbol{\varphi}_{j}d\Omega &= \int_{\Omega} \mathbf{q}^{n} \cdot \boldsymbol{\varphi}_{j} \mathrm{d}\Omega \nonumber \\
	&- a_{21}\Delta t \sum_{K \in \mathcal{T}_{\mathcal{H}}}\int_{\partial K}\mathbf{F}\rpth{\mathbf{q}^{n}}\bm{n} \cdot \boldsymbol{\varphi}_{j} \mathrm{d}\Sigma - a_{21}\Delta t\int_{\Omega}\mathbf{B}\rpth{\mathbf{q}^{n}}\grad\mathbf{q}^{n} \cdot \boldsymbol{\varphi}_{j} \mathrm{d}\Omega \nonumber \\
	&+ \tilde{a}_{21}\Delta t \int_{\Omega} \mathbf{S}_{\mathrm{I}}\rpth{\mathbf{q}^{n}} \cdot \boldsymbol{\varphi}_{j} \mathrm{d}\Omega + \tilde{a}_{22}\Delta t \int_{\Omega} \mathbf{S}_{\mathrm{I}}\rpth{\mathbf{q}^{(n,2)}} \cdot \boldsymbol{\varphi}_{j} \mathrm{d}\Omega \\
	&+ a_{21}\Delta t \int_{\Omega} \mathbf{S}_{\mathrm{E}} \cdot \boldsymbol{\varphi}_{j} \mathrm{d}\Omega. \nonumber
\end{align}
Here $\boldsymbol{\varphi}_{j}, j = 1, \dots, \text{dim}(\bm{Q}_{0})$ denotes the set of basis functions of the space $\bm{Q}_{0}$, with $\left|\mathcal{T}_{H}\right|$ being the number of elements of the computational mesh. Moreover, $\bm{n}$ is the outward unit normal vector to the boundary $\partial K$ of the element $K$. Similarly, the discrete weak form associated to the third stage of the IMEX reads as follows
\begin{align}\label{eq:eq_stage_3_up}
	\int_{\Omega} \mathbf{q}^{(n,3)} \cdot \tilde{\boldsymbol{\varphi}}_{j} \mathrm{d}\Omega &= \int_{\Omega} \mathbf{q}^{n} \cdot \tilde{\boldsymbol{\varphi}}_{j} \mathrm{d}\Omega \nonumber \\
	&+ a_{31}\Delta t \int_{\Omega} \mathbf{F}\rpth{\mathbf{q}^{n}} : \grad\tilde{\boldsymbol{\varphi}}_{j} \mathrm{d}\Omega - a_{31}\Delta t \int_{\partial\Omega} \mathbf{F}\rpth{\mathbf{q}^{n}}\bm{n}_{\Omega} \cdot \tilde{\boldsymbol{\varphi}}_{j} \mathrm{d}\Sigma \nonumber \\
	&- a_{31}\Delta t \int_{\Omega} \mathbf{B}\rpth{\mathbf{q}^{n}}\grad\mathbf{q}^{n} \cdot \tilde{\boldsymbol{\varphi}}_{j}\mathrm{d}\Omega \nonumber \\
	&+ a_{32}\Delta t \int_{\Omega} \mathbf{F}\rpth{\mathbf{q}^{(n,2)}} : \grad\tilde{\boldsymbol{\varphi}}_{j} \mathrm{d}\Omega - a_{32}\Delta t \int_{\partial\Omega} \mathbf{F}\rpth{\mathbf{q}^{(n,2)}}\bm{n}_{\Omega} \cdot \tilde{\boldsymbol{\varphi}}_{j} \mathrm{d}\Sigma \\
	&- a_{32}\Delta t \int_{\Omega} \mathbf{B}\rpth{\mathbf{q}^{(n,2)}}\grad\mathbf{q}^{(n,2)} \cdot \tilde{\boldsymbol{\varphi}}_{j} \mathrm{d}\Omega \nonumber \\
	&+ \tilde{a}_{31}\Delta t \int_{\Omega} \mathbf{S}_{\mathrm{I}}\rpth{\mathbf{q}^{n}} \cdot \tilde{\boldsymbol{\varphi}}_{j} \mathrm{d}\Omega + \tilde{a}_{32}\Delta t \int_{\Omega} \mathbf{S}_{\mathrm{I}}\rpth{\mathbf{q}^{(n,2)}} \cdot \tilde{\boldsymbol{\varphi}}_{j} \mathrm{d}\Omega \nonumber \\
	&+ \gamma\Delta t \int_{\Omega} \mathbf{S}_{\mathrm{I}}\rpth{\mathbf{q}^{(n,3)}} \cdot \tilde{\boldsymbol{\varphi}}_{j} \mathrm{d}\Omega \nonumber \\
	&+ a_{31}\Delta t \int_{\Omega} \mathbf{S}_{\mathrm{E}} \cdot \tilde{\boldsymbol{\varphi}}_{j} \mathrm{d}\Omega + a_{32}\Delta t \int_{\Omega} \mathbf{S}_{\mathrm{E}} \cdot \tilde{\boldsymbol{\varphi}}_{j} \mathrm{d}\Omega. \nonumber
\end{align}
Here, $\tilde{\boldsymbol{\varphi}}_{i}, j = 1, \dots, \text{dim}(\tilde{\bm{Q}}_{1})$ denotes the set of basis functions of the space $\tilde{\bm{Q}}_{1}$, while $\bm{n}_{\Omega}$ is the outward unit normal vector to the boundary $\partial\Omega$. Notice that $\grad\mathbf{q}^{(n,2)}$ is not defined in the space $\bm{Q}_{0}$. In Section~\ref{ssec:non_conservative_terms} we will describe the strategy adopted for the discretization of the non-conservative term in this case. Finally, the discrete weak form associated with the final update of the IMEX method can be expressed in a similar manner.

%%%%%%%%%% Discretization of non-conservative terms %%%%%%%
\subsection{Discretization of non-conservative terms}
\label{ssec:non_conservative_terms}

In this section, we analyze the treatment of the non-conservative term of system~\eqref{eq:lava_model_Z} for discontinuous spaces. Following the work of~\cite{pares:2006}, which is based on the so-called Dal Maso, Le Floch and Murat (DLM) theory~\cite{dalmaso:1995}, the Path-Conservative (PC) formulation reads as follows
\begin{equation}
	\int_{K} \mathbf{B}\rpth{\mathbf{q}}\hat{\nabla}\mathbf{q} \cdot \tilde{\boldsymbol{\varphi}} \mathrm{d}\Omega = \int_{K} \mathbf{B}\rpth{\mathbf{q}}\nabla\mathbf{q} \cdot \tilde{\boldsymbol{\varphi}} \mathrm{d}\Omega + \int_{\partial K} \boldsymbol{\mathcal{D}}\rpth{\mathbf{q}^{-}, \mathbf{q}^{+}, \bm{n}} \cdot \tilde{\boldsymbol{\varphi}} \mathrm{d}\Sigma.
\end{equation}
In the particular case $\mathbf{q} \in \bm{Q}_{0}$, one obtains
\begin{equation}
	\int_{K} \mathbf{B}\rpth{\mathbf{q}}\hat{\nabla}\mathbf{q} \cdot \tilde{\boldsymbol{\varphi}}d\Omega = \int_{\partial K} \boldsymbol{\mathcal{D}}\rpth{\mathbf{q}^{-}, \mathbf{q}^{+}, \bm{n}} \cdot \tilde{\boldsymbol{\varphi}} \mathrm{d}\Sigma.
\end{equation}
The function $\boldsymbol{\mathcal{D}}$ depends on the value of $\mathbf{q}$ from the interior $\mathbf{q}^{-}$, the value of $\mathbf{q}$ from the exterior $\mathbf{q}^{+}$, and on the outward unit normal vector to the element $K$. The function $\boldsymbol{\mathcal{D}}$ is such that
\begin{subequations}
	\begin{align}
		&\boldsymbol{\mathcal{D}}\rpth{\mathbf{q}, \mathbf{q}, \bm{n}} = 0 \quad \forall \mathbf{q}, \bm{n}, \label{eq:PC_consistency} \\
		&\boldsymbol{\mathcal{D}}\rpth{\mathbf{q}^{-}, \mathbf{q}^{+}, \bm{n}} + \boldsymbol{\mathcal{D}}\rpth{\mathbf{q}^{+}, \mathbf{q}^{-}, -\bm{n}} = \spth{\int_{0}^{1} \mathbf{B}\rpth{\boldsymbol{\Psi}\rpth{\mathbf{q}^{-}, \mathbf{q}^{+}, s}}\pad{\boldsymbol{\Psi}}{s} \mathrm{d}s}\bm{n}, \label{eq:PC_conservative}
	\end{align}
\end{subequations}
The function $\boldsymbol{\Psi}$ denotes a family of paths connecting $\mathbf{q}^{-}$ and $\mathbf{q}^{+}$ such that 
\begin{equation*}
	\boldsymbol{\Psi}\rpth{\mathbf{q}^{-}, \mathbf{q}^{+}, 0} = \mathbf{q}^{-}, \qquad \boldsymbol{\Psi}\rpth{\mathbf{q}^{-}, \mathbf{q}^{+}, 1} = \mathbf{q}^{+}, \qquad \boldsymbol{\Psi}\rpth{\mathbf{q}, \mathbf{q}, s} = \mathbf{q}.
\end{equation*}
We refer to~\cite{pares:2006} for a detailed discussion of the regularity required by $\boldsymbol{\Psi}$. The standard choice for the path is the segment~\cite{busto:2021, gaburro:2024, gatti:2025}, i.e.
\begin{equation}\label{eq:linear_paths}
	\boldsymbol{\Psi}\rpth{\mathbf{q}^{-}, \mathbf{q}^{+}, s} = \mathbf{q}^{-} + s\rpth{\mathbf{q}^{+} - \mathbf{q}^{-}}.
\end{equation}
Following~\cite{busto:2021, gatti:2025}, given \eqref{eq:non_cons_B}, we consider the choice
\begin{equation}\label{eq:PC_trapezoidal}
	\boldsymbol{\mathcal{D}}\rpth{Z^{-}, Z^{+}, \bm{n}} = \frac{1}{2}\spth{\int_{0}^{1} gh\rpth{\boldsymbol{\Psi}\rpth{\mathbf{q}^{-}, \mathbf{q}^{+}, s}}\pad{\boldsymbol{\Psi}}{s} \mathrm{d}s}\bm{n}
\end{equation}
and we approximate this integral by means of a trapezoidal rule. In particular, if $\mathbf{q} \in \bm{Q}_{0}$ as for $\mathbf{q}^{(n,2)}$, we obtain
\begin{align}
	\sum_{K} \int_{K} gh\hat{\nabla}Z \cdot \tilde{\boldsymbol{\varphi}}d\Omega &= \sum_{e \in \mathcal{E}} \rpth{Z^{+} - Z^{-}} \int_{e}\spth{\int_{0}^{1} g\rpth{h^{-} + s\rpth{h^{+} - h^{-}}}ds}\bm{n} \cdot \tilde{\boldsymbol{\varphi}} \mathrm{d}\Sigma \nonumber \\
	&\approx \sum_{e \in \mathcal{E}} \rpth{Z^{+} - Z^{-}}\frac{1}{2}g\rpth{h^{-} + h^{+}} \int_{e}\tilde{\boldsymbol{\varphi}} \cdot \bm{n} \mathrm{d}\Sigma.
\end{align}

Finally, we emphasize that preserving the well-balancing property is essential both in the computation of the intermediate stage and in the solution update. Indeed, a discretization of the non-conservative terms that is not fully consistent with the approximation of the conservative fluxes would generate an artificial imbalance of forces, thereby destroying the well-balancing property. Implementation details, as well as a proof of the discrete well-balancing property for well-prepared data, are provided in Appendix \ref{app:detail_wb}.

%%%%%%%%%%% von Neumann stability analysis %%%%%%%%%%%
\subsection{Completion of the IMEX scheme for the (pseudo-)staggered Galerkin method}
\label{ssec:staggered_von_Neumann_et_al}

The final step to complete the (pseudo-)staggered Galerkin method is to determine the remaining coefficients of the IMEX scheme (see Table~\ref{tab:imex_rk2_butch}). As already mentioned in Section~\ref{ssec:novel_IMEX}, the starting point to achieve this goal is to perform a von Neumann stability analysis~\cite{deriaz:2020}. For this purpose, we consider the one-dimensional linear advection-reaction problem with periodic boundary conditions
\begin{equation}\label{eq:test_equation_von_Neumann}
	\pad{q}{t} + a\pad{q}{x} + \chi q = 0, \qquad x \in \rpth{0, L}, \quad t > 0.
\end{equation}
Here, $q$ is a scalar variable, while $\chi$ and $a$ are real constant coefficients. We consider a uniform mesh of size $\Delta x$. The classical mass lumping technique~\cite{hughes:2003, zienkiewicz:2005} is adopted when dealing with solutions in continuous space $\tilde{Q}_{1}$, to avoid building a global matrix. The resulting fully discrete stages read therefore as follows
\begin{subequations}\label{eq:fully_discrete_von_Neumann}
	\begin{align}
		q_{j+\frac{1}{2}}^{(n,2)} &= \frac{q_{j}^{n} + q_{j+1}^{n}}{2} - \tilde{a}_{21} \Phi\frac{q_{j}^{n} + q_{j+1}^{n}}{2} - \gamma\Phi q_{j+\frac{1}{2}}^{(n,2)} - a_{21}\nu\rpth{q_{j+1}^{n} - q_{j}^{n}}, \label{eq:fully_discrete_von_Neumann_stage2} \\
		q_{j}^{(n,3)} &= q_{j}^{n} - \tilde{a}_{31}\Phi q_{j}^{n} -\tilde{a}_{32}\Phi\frac{q_{j-\frac{1}{2}}^{(n,2)} + q_{j+\frac{1}{2}}^{(n,2)}}{2} - \gamma\Phi q_{j}^{(n,3)} \nonumber \\
		&- a_{31}\nu\dfrac{q_{j+1}^{n} - q_{j-1}^{n}}{2} + a_{32}\nu\rpth{q_{j-\frac{1}{2}}^{(n,2)} - q_{j+\frac{1}{2}}^{(n,2)}}, \label{eq:fully_discrete_von_Neumann_stage3} \\
		q_{j}^{n+1} &= q_{j}^{n} - \tilde{a}_{31}\rpth{\Phi q_{j}^{n} + \nu\frac{q_{j+1}^{n} - q_{j-1}^{n}}{2}} \nonumber \\ &-\tilde{a}_{32}\rpth{\Phi\frac{q_{j-\frac{1}{2}}^{(n,2)} + q_{j+\frac{1}{2}}^{(n,2)}}{2} - \nu\rpth{q_{j-\frac{1}{2}}^{(n,2)} - q_{j+\frac{1}{2}}^{(n,2)}}} \nonumber \\
		&- \gamma\rpth{\Phi q_{j}^{(n,3)} + \nu\frac{q_{j+1}^{(n,3)} - q_{j-1}^{(n,3)}}{2}}, \label{eq:fully_discrete_von_Neumann_update}
	\end{align}
\end{subequations}
with
$$\nu = \frac{a\Delta t}{\Delta x}, \qquad \Phi = \chi\Delta t.$$
In order to perform the classical von Neumann stability analysis, one assumes that $q_{j}^{0} = \exp\rpth{ik x_{j}}$, where $i$ is the imaginary unit, $k$ represents the wave number and $x_{j} = j\Delta x$. Before proceeding, we provide the following
%.
\begin{definition}\label{def:space_time_L_stability}
	A numerical method for~\eqref{eq:test_equation_von_Neumann} is said to be space-time \texttt{L}-stable if
	\begin{equation}
		\lim_{\Phi \to \infty} G\rpth{\nu, \Phi, \theta} = 0,
	\end{equation}
	where $\theta = k\Delta x$ and $G$ denotes the amplification factor of a von Neumann stability analysis.
\end{definition}
Notice that, to the best of our knowledge, this represents a novel definition of stability and is one of the main contributions of the present work. After computing $q_{j}^{1}$ with~\eqref{eq:fully_discrete_von_Neumann}, tedious calculations (see Appendix~\ref{app:staggered_von_Neumann_computations} for the main steps) show that the amplification factor in the limit case $\Phi \to \infty$ reduces to
\begin{equation}\label{eq:Glim}
	G_{\mathrm{lim}} = \frac{1}{2 \gamma^{2}} \spth{i \gamma \nu \sin{\theta}\rpth{\tilde{a}_{21}\rpth{\tilde{a}_{32} - 2a_{32}} - \tilde{a}_{21}\tilde{a}_{32}\cos{\theta} + 2a_{31}\gamma} + \tilde{a}_{21}\tilde{a}_{32}\cos{\theta} + \tilde{a}_{21}\tilde{a}_{32} - 2\tilde{a}_{31}\gamma}.
\end{equation}

To avoid any dependence of $G_{\mathrm{lim}}$ on the phase $\theta$, as well as to satisfy the space-time \texttt{L}-stability property introduced in Definition~\eqref{def:space_time_L_stability}, the following relations have to be satisfied
\begin{equation}\label{eq:IMEX_conditions_von_Neumann}
	\begin{cases}
		\tilde{a}_{21}\tilde{a}_{32}= 0, \\
		\tilde{a}_{21}\rpth{\tilde{a}_{32} - 2a_{32}}+ 2a_{31}\gamma = 0, \\
		\tilde{a}_{31}\gamma = 0.
	\end{cases}
\end{equation}
To derive our method, we proceed by considering the space-time \texttt{L}-stability conditions listed above together with the conditions reported in Equation~\eqref{eq:number_1_cond} and~\eqref{eq:number_2_cond}. The third equation in~\eqref{eq:IMEX_conditions_von_Neumann} immediately yields $\tilde{a}_{31} = 0$. From the first equation listed in~\eqref{eq:IMEX_conditions_von_Neumann}, we can choose either $\tilde{a}_{21} = 0$ (i) or $\tilde{a}_{32} = 0$ (ii). In the case (i), we obtain $a_{31} = 0$, $\gamma = 1 - \sqrt{2}/2$, $\tilde{a}_{32} = \sqrt{2}/2$, and \eqref{eq:number_2_cond} reduces to
\begin{equation}\label{eq:supplmentary_condition_IMEX_von_Neumann}
	-a_{21}\sqrt{2} - a_{32}\rpth{2 - \sqrt{2}} + 1 = 0.
\end{equation}
On the other hand, one can easily verify that case (ii) does not lead anywhere, since it results in a contradiction. Indeed, imposing $\tilde{a}_{31} = \tilde{a}_{32} = 0$ implies $\gamma = 1$ because of~\eqref{eq:order1_condition}, but this immediately contradicts~\eqref{eq:a32_tilde} for which $\tilde{a}_{32} = 0$ if and only if $\gamma = \frac{1}{2}$. The coefficients of the resulting implicit method are therefore
\begin{table}[H]
	\centering
	\begin{tabular}{c|ccc}
		$0$ & $0$ & $0$ & $0$ \\
		$1 - \frac{\sqrt{2}}{2}$ & $0$ & $1 - \frac{\sqrt{2}}{2}$ & $0$ \\
		$1$ & $0$ & $\frac{\sqrt{2}}{2}$ & $1 - \frac{\sqrt{2}}{2}$ \\[2pt]
		\hline
		\rule{0pt}{3ex} & $0$ & $\frac{\sqrt{2}}{2}$ & $1 - \frac{\sqrt{2}}{2}$
	\end{tabular}
\end{table}
The implicit scheme is of type ARS~\cite{ascher:1997} and has already been proposed as a two-stage SDIRK method in~\cite{alexander:1977, crouzeix:1976} and used in~\cite{luo:2005} for unconstrained optimization. It is also \texttt{L}-stable in the classical sense related to the discretization of ODEs~\cite{kennedy:2016}.

We are left with one free parameter for the explicit companion scheme to complete the derivation of the IMEX method. One could impose $\mathbf{c} = \tilde{\mathbf{c}}$, from which we obtain
\begin{equation}
	a_{21} = \gamma = 1 - \frac{\sqrt{2}}{2}, \qquad a_{32} = 1.
\end{equation}
It is immediately clear that these two values satisfy~\eqref{eq:supplmentary_condition_IMEX_von_Neumann}. Hence, the resulting IMEX scheme has Butcher tableaux reported in Table~\ref{tab:IMEX_c_eq_ctilde}.
\begin{table}[H]
	\begin{minipage}{0.45\textwidth}
		\begin{center}
			\begin{tabular}{c|ccc}
				$0$ & $0$ & $0$ & $0$ \\
				$1-\frac{\sqrt{2}}{2}$ & $1-\frac{\sqrt{2}}{2}$ & $0$ & $0$ \\
				$1$ & $0$ & $1$ & $0$ \\[2pt]
				\hline
				\rule{0pt}{3ex} & $0$ & $\frac{\sqrt{2}}{2}$ & $1-\frac{\sqrt{2}}{2}$
			\end{tabular}
		\end{center}
	\end{minipage}
	\begin{minipage}{0.45\textwidth}
		\begin{center}
			\begin{tabular}{c|ccc}
				$0$ & $0$ & $0$ & $0$ \\
				$1-\frac{\sqrt{2}}{2}$ & $0$ & $1-\frac{\sqrt{2}}{2}$ & $0$ \\
				$1$ & $0$ & $\frac{\sqrt{2}}{2}$ & $1-\frac{\sqrt{2}}{2}$ \\[2pt]
				\hline
				\rule{0pt}{3ex} & $0$ & $\frac{\sqrt{2}}{2}$ & $1-\frac{\sqrt{2}}{2}$
			\end{tabular}
		\end{center}
	\end{minipage}
	\caption{Butcher tableaux of the the space-time \texttt{L}-stable IMEX scheme with $\mathbf{c} = \tilde{\mathbf{c}}$. Left: explicit method. Right: implicit method.}
	\label{tab:IMEX_c_eq_ctilde}
\end{table}
\noindent
Imposing $\mathbf{c} = \tilde{\mathbf{c}}$ can be useful because it closes the system directly without the need to impose any other condition. However, we note that this is not strictly mandatory. Hence, we choose to determine the remaining free coefficient of the explicit companion scheme driven by the leading-order term of the diffusion error. For this purpose, we compute the modified equation of the fully discrete scheme by substituting the Taylor series expansions in space and time, so that we can isolate the local truncation error of the scheme. The time derivatives are written as functions of spatial derivatives with the help of the Lax--Wendroff procedure, also known as the Cauchy--Kovalevskaya procedure for higher-order methods~\cite{harten:1997}. We obtain the following modified equation where we discard the contributions greater than the second order
\begin{align}\label{eq:modified_equation}
	&\pad{q}{t} + a\pad{q}{x} + \chi q = \rpth{\frac{\sqrt{2}}{2} - \frac{2}{3}}\Delta t^{2}\chi^{3}q \\
	&+ \spth{\rpth{2 - \frac{3}{\sqrt{2}}}a_{32}^2 - \rpth{\frac{1}{2} - \frac{\sqrt{2}}{2}}a_{32} + \rpth{\frac{1}{2} - \frac{\sqrt{2}}{2} + \frac{1}{4\sqrt{2}\nu^{2}}}}a^{2}\Delta t^{2} \chi \frac{\partial^{2}q}{\partial x^{2}} \nonumber \\
	&+\spth{\rpth{2 - \frac{3}{\sqrt{2}}}a_{32}^{2} - \rpth{\frac{\frac{1}{2} - \sqrt{2}}{2}} a_{32} + \frac{1}{6}\rpth{1 - \frac{1}{\nu^{2}}} }a^{3}\Delta t^{2}\frac{\partial^{3}q}{\partial x^{3}}. \nonumber
\end{align}
The reader can refer to the Mathematica~\cite{Mathematica2024} notebook \texttt{LW\_proc.nb} in the GitHub repository \url{https://github.com/federicg/lava-flow.git} for the actual computations. 

First of all, we point out that the scheme is of second order. Indeed, the modified equation shows that the numerical scheme solves the original advection-reaction equation up to error terms that are all proportional to $\Delta t^{2}$, multiplying spatial derivatives of $q$ of order two and three through the factors $a^{2}$ and $a^{3}$, respectively. Since the Courant number $\nu = a\Delta t/\Delta x$ is kept fixed (or bounded) in the present analysis, the leading-order error terms in Equation~\eqref{eq:modified_equation} are equivalently of order $\mathcal{O}(\Delta t^{2}) = \mathcal{O}(\Delta x^{2})$, confirming that the truncation error of the scheme is second order in both time and space, i.e., the scheme is second-order accurate in space-time. This result is achieved in spite of the use of the (pseudo-)staggered method for the space discretization. 

Next, by imposing the stability condition on the diffusion coefficient of the modified equation, i.e.,
\begin{equation}\label{eq:diffusion_modified_equation}
	\spth{\rpth{2 - \frac{3}{\sqrt{2}}}a_{32}^2 - \rpth{\frac{1}{2} - \frac{\sqrt{2}}{2}}a_{32} + \rpth{\frac{1}{2} - \frac{\sqrt{2}}{2} + \frac{1}{4\sqrt{2}\nu^{2}}}}a^{2}\Delta t^{2} \chi \ge 0, 
\end{equation}
we obtain, for $\chi > 0$,
\begin{equation}
	\left|\nu\right| \le \frac{\sqrt{\frac{1}{3 \sqrt{2} a_{32}^2-4 a_{32}^2-\sqrt{2} a_{32} + a_{32} + \sqrt{2}-1}}}{2^{3/4}}.
\end{equation}
The coefficient $a_{32}$ is then set by maximizing the value of $|\nu|$. This, together with \eqref{eq:supplmentary_condition_IMEX_von_Neumann}, leads to
$$a_{32} = \frac{\sqrt{2} - 1}{2\rpth{3\sqrt{2} - 4}}, \qquad a_{21} = \frac{\sqrt{2}}{4}.$$
Then, we have the following bound on $|\nu|$
\begin{equation}
	\left|\nu\right| \le \sqrt{\frac{2\rpth{17\sqrt{2} - 24}}{215\sqrt{2} - 304}} \approx 1.22.
\end{equation}
Finally, the proposed second-order space-time \texttt{L}-stable IMEX scheme has coefficients reported in Table~\ref{tab:IMEX_maximized_stability}.

\begin{figure}[h!]
	\centering
	\includegraphics[scale = 0.5]{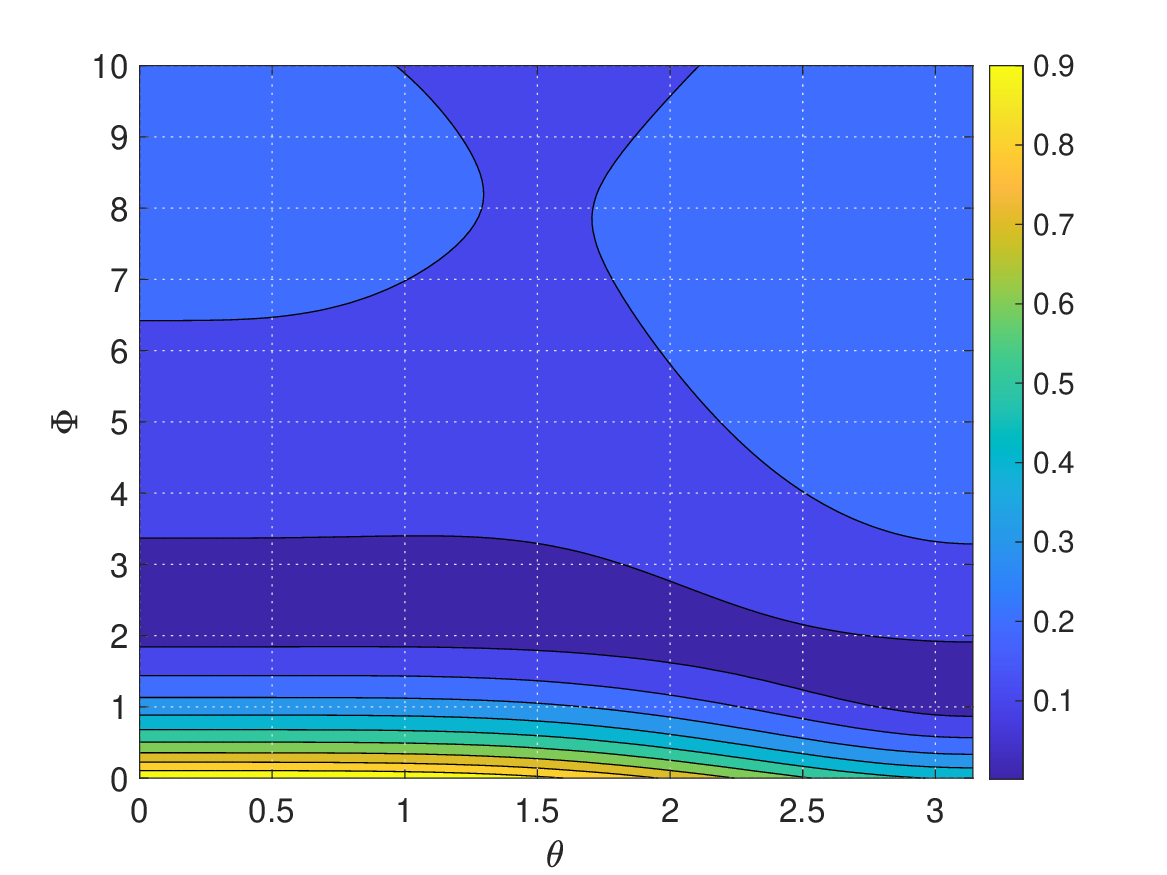}
	\caption{Contour plot of $\left|G_{\mathrm{lim}}\right|$~\eqref{eq:Glim} for $|\nu|\approx 1.22$.}
	 \label{fig:contourPlot_Phi_theta}
\end{figure}

\begin{table}[H]
	\begin{minipage}{0.475\textwidth}
		\begin{center}
			\begin{tabular}{c|ccc}
				$0$ & $0$ & $0$ & $0$ \\
				$\frac{\sqrt{2}}{4}$ & $\frac{\sqrt{2}}{4}$ & $0$ & $0$ \\
				$\frac{\sqrt{2} - 1}{2\rpth{3\sqrt{2} - 4}}$ & $0$ & $\frac{\sqrt{2} - 1}{2\rpth{3\sqrt{2} - 4}}$ & $0$ \\
				[2pt]\hline
				\rule{0pt}{3ex} & $0$ & $\frac{\sqrt{2}}{2}$ & $1 - \frac{\sqrt{2}}{2}$
			\end{tabular}
		\end{center}
	\end{minipage}
	\begin{minipage}{0.475\textwidth}
		\begin{center}
			\begin{tabular}{c|ccc}
				$0$ & $0$ & $0$ & $0$ \\
				$1 - \frac{\sqrt{2}}{2}$ & $0$ & $1 - \frac{\sqrt{2}}{2}$ & $0$ \\
				$1$ & $0$ & $\frac{\sqrt{2}}{2}$ & $1 - \frac{\sqrt{2}}{2}$ \\
				[2pt]\hline
				\rule{0pt}{3ex} & $0$ & $\frac{\sqrt{2}}{2}$ & $1 - \frac{\sqrt{2}}{2}$
			\end{tabular}
		\end{center}
	\end{minipage}
	\caption{Butcher tableaux of the the space-time \texttt{L}-stable IMEX scheme obtained maximizing $\left|\nu\right|$. Left: explicit method. Right: implicit method.}
	\label{tab:IMEX_maximized_stability}
\end{table}
\noindent
No restrictions are experienced from a stability point of view for $\Phi \ge 0 $ as $\left|G_{\mathrm{lim}}\right| \le 1$ for all values of $\theta$ (Figure~\ref{fig:contourPlot_Phi_theta}). In particular, one can easily verify that $\left|G_{\mathrm{lim}}\right| = 1$ for a pure advection problem, i.e., for $\Phi = 0$ in~\eqref{eq:test_equation_von_Neumann}.

It is worth noting that the schemes reported in Table~\ref{tab:IMEX_c_eq_ctilde} and Table~\ref{tab:IMEX_maximized_stability} share the same implicit companion method and differ only in the explicit scheme, which in the first case is determined by imposing $\mathbf{c} = \tilde{\mathbf{c}}$, whereas in the second case it is obtained by maximizing the stability according to the Lax--Wendroff procedure.

While extending the present analysis to higher-order methods lies beyond the scope of this work, the use of IMEX schemes naturally raises the question of whether a similar approach can be pursued in that setting. The main difficulty would likely be related to the larger number of free parameters involved; consequently, although we have not verified this, the outcome of the von Neumann stability analysis may not be sufficient to fully determine the implicit companion method, as is the case for the three-stage second-order schemes considered here. We briefly discuss this issue in the concluding remarks of Section~\ref{sec:conclu}.

%%%%%%%%%%%%%% Numerical results %%%%%%%%%%%%%%%
\section{Numerical results}
\label{sec:num}

We consider a number of relevant test cases to assess the properties of the numerical method depicted in Sections~\ref{sec:time_disc} and~\ref{sec:space_disc}. 
In the following, to choose the time step for the numerical simulations of system~\eqref{eq:lava_model_Z}, we consider the maximum Courant number along each direction.
Following~\cite{gatti:2023a, gatti:2025, gatti:2023b}, we set a numerical threshold $h_{\mathrm{min}}$ on the height of the material $h$, to provide a discretization of the wet region $\Omega_{w}$. More specifically, in the regions where $h$ is smaller than $h_{\mathrm{min}}$, the model equations are modified are as follows
\begin{equation}
	\pad{h}{t} = 0, \quad h\vel = \bm{0}, \quad hT = 0.
\end{equation}
Notice that this guaranties conservation of mass since we act only on the model fluxes and do not modify the material height in case it is lower than $h_{\mathrm{min}}$. We set $h_{\mathrm{min}} = \SI[parse-numbers=false]{10^{-5}}{\meter}$ for the simulations in which a wet-dry interface is present. Finally, following~\cite{gatti:2025}, we approximate the Dirac's delta appearing in the source term through the following Gaussian function
\begin{equation}\label{eq:Gaussian_lava}
	f_{v}(r) = \dfrac{1}{2\pi\sigma}e^{-r^{2}/(2\sigma)}, \quad \text{with } r^{2} = (x - x_{v})^{2} + (y - y_{v})^{2},
\end{equation}
where $\sigma$ is a parameter that has the dimension of an area and is defined by the user. We recall that due to the nature of the Gaussian function, we have decided to compute exactly the integrals resulting from the projection $L^{2}$ of the Gaussian function on the adopted finite element spaces. In particular, the projection onto the set of bilinear basis functions requires the use of the integration by parts, leading to the presence of the error function. We employ the error function implemented in the C++ standard library.

\subsection{One-dimensional tests}
\label{ssec:1D_tests}

In this section, we consider the linear advection-reaction equation~\eqref{eq:test_equation_von_Neumann} to provide a first validation of the proposed nu\-me\-ric\-al me\-thod. First, we focus on a model problem characterized by the sole presence of the reaction term, while in a second test, we also include the advection term, and we provide a convergence analysis on smooth solutions.

\subsubsection{Linear reaction problem}
\label{ssec:linear_reaction}

\begin{figure}[h!] 
	\centering  
	\includegraphics[scale=0.32]{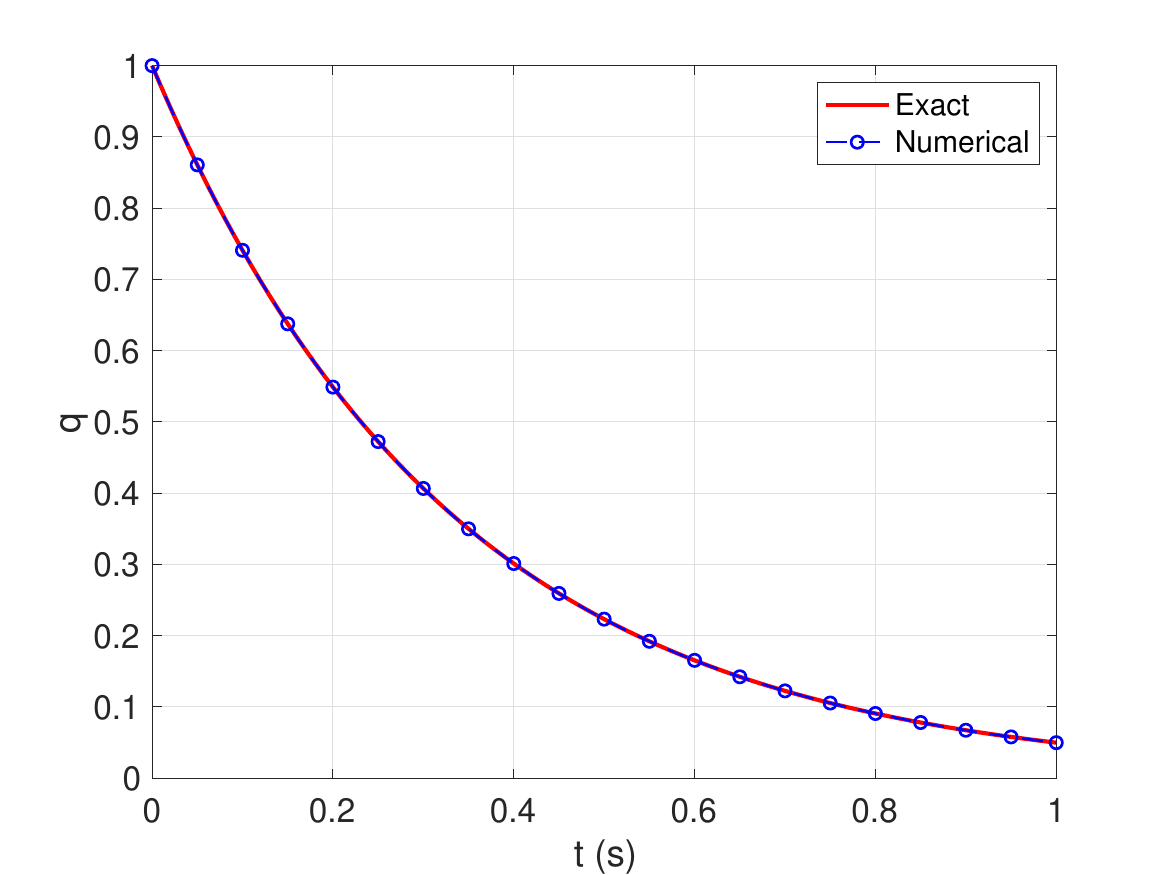} 
	\includegraphics[scale=0.32]{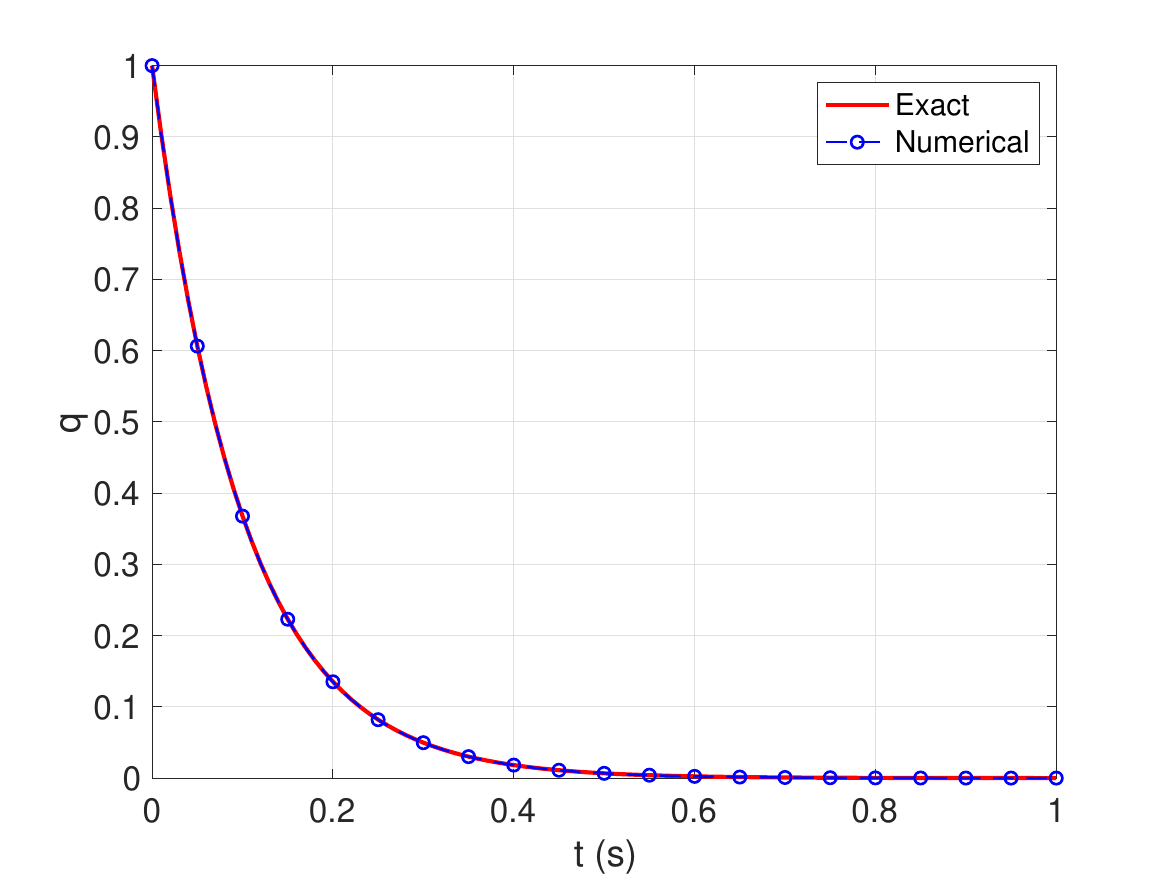}
	\includegraphics[scale=0.32]{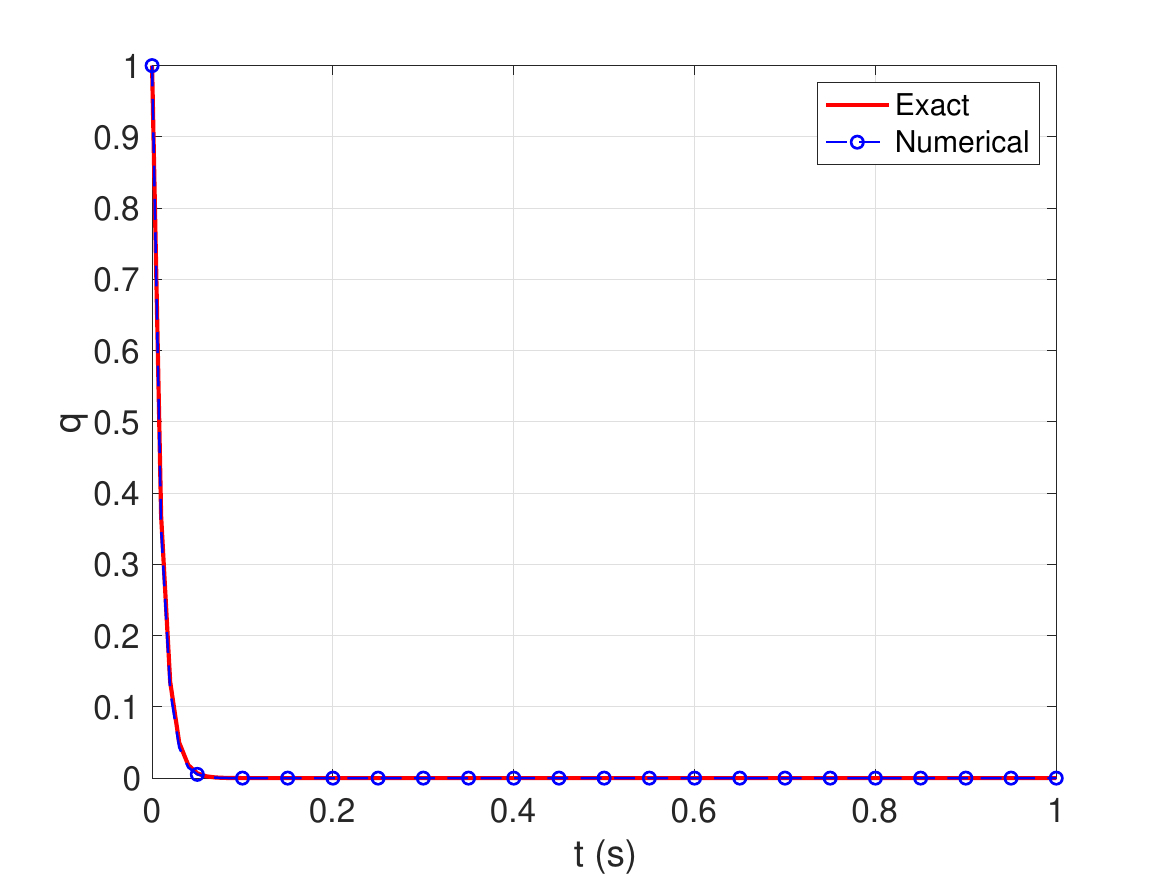}
	\includegraphics[scale=0.32]{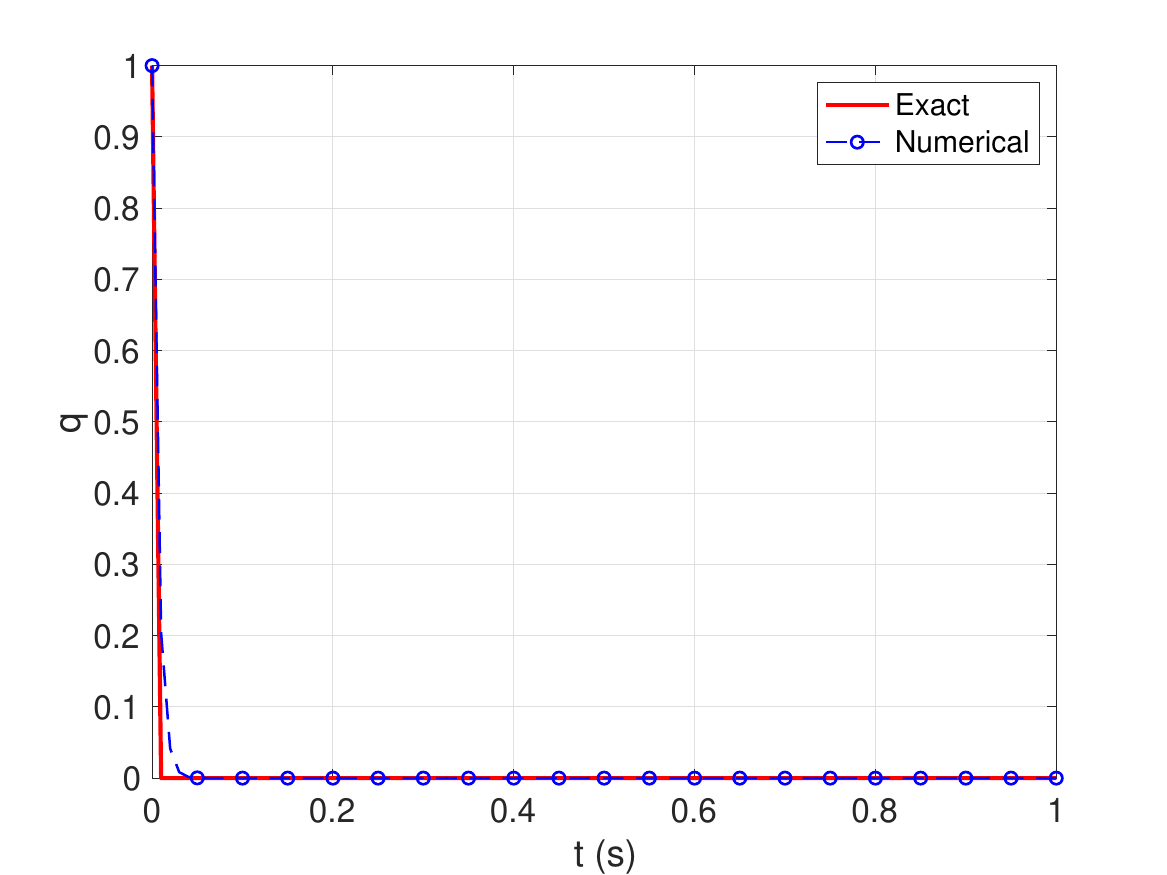}
	\caption{Linear reaction test case, comparison between the exact solution (solid red lines) and numerical results (blue dots). Top-left: $\chi = 3$. Top-right: $\chi = 10$. Bottom-left: $\chi = 100$. Bottom-right: $\chi = 1000$.}
	\label{fig:linear_reaction}
\end{figure}

The present test case has been considered, e.g., by~\cite{dumbser:2008} as a benchmark for numerical methods for stiff problems. This test considers null advection velocity, so that~\eqref{eq:test_equation_von_Neumann} reduces to
\begin{equation}\label{eq:lin_reaction_eq}
	\frac{\partial q}{\partial t} + \chi q = 0, \quad t > 0.
\end{equation}
The analytical solution of~\eqref{eq:lin_reaction_eq} is 
\begin{equation}
	q(t; \chi) = q_{0}e^{-\chi t},
\end{equation}
with $q_{0}$ denoting the initial condition. The problem becomes more stiff as the scalar value $\chi$ increases and, as $\chi \to \infty$, the solution becomes discontinuous. Following~\cite{dumbser:2008}, we set $q_{0} = 1$ in the entire domain, while a fixed time step $\Delta t = \SI{0.01}{\second}$ is used and a total of uniform cells $300$ are used. We verify that the numerical method developed reproduces the behavior of the analytical solution for $\chi \in \cpth{3, 10, 100, 1000}$. One can easily notice that an excellent agreement is established with the analytical solution for all values in the range of $\chi$ (Figure~\ref{fig:linear_reaction}). We recall that no TVD technique is employed to avoid spurious oscillations that arise when high-order methods are employed in the case of discontinuous solutions. In fact, despite the absence of the convection term, the numerical scheme still has a spatial dependence given by the discretization of the reaction term (see Equation~\eqref{eq:fully_discrete_von_Neumann}).

\subsubsection{Linear advection-reaction problem}
\label{ssec:linear_advection_reaction}

\begin{figure}[h!]
	\centering
	\includegraphics[scale = 0.32]{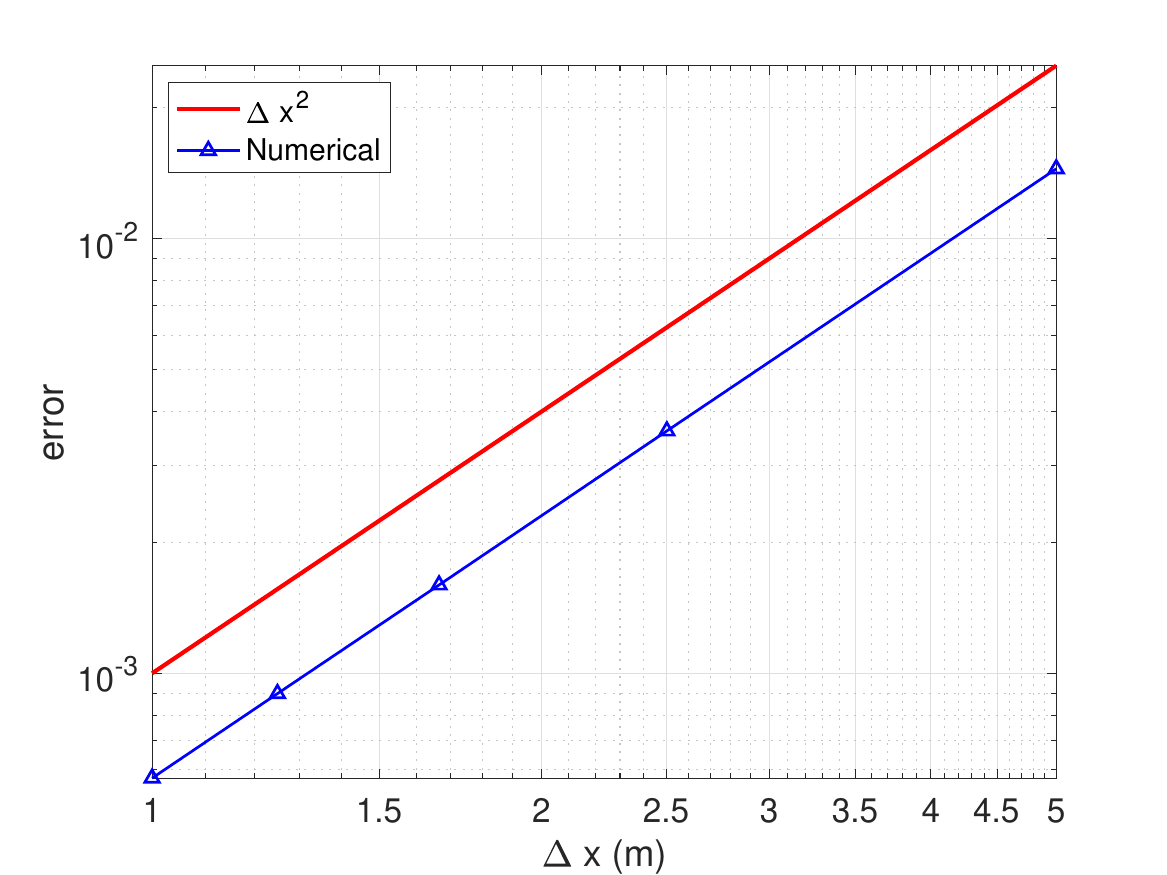}
	\includegraphics[scale = 0.32]{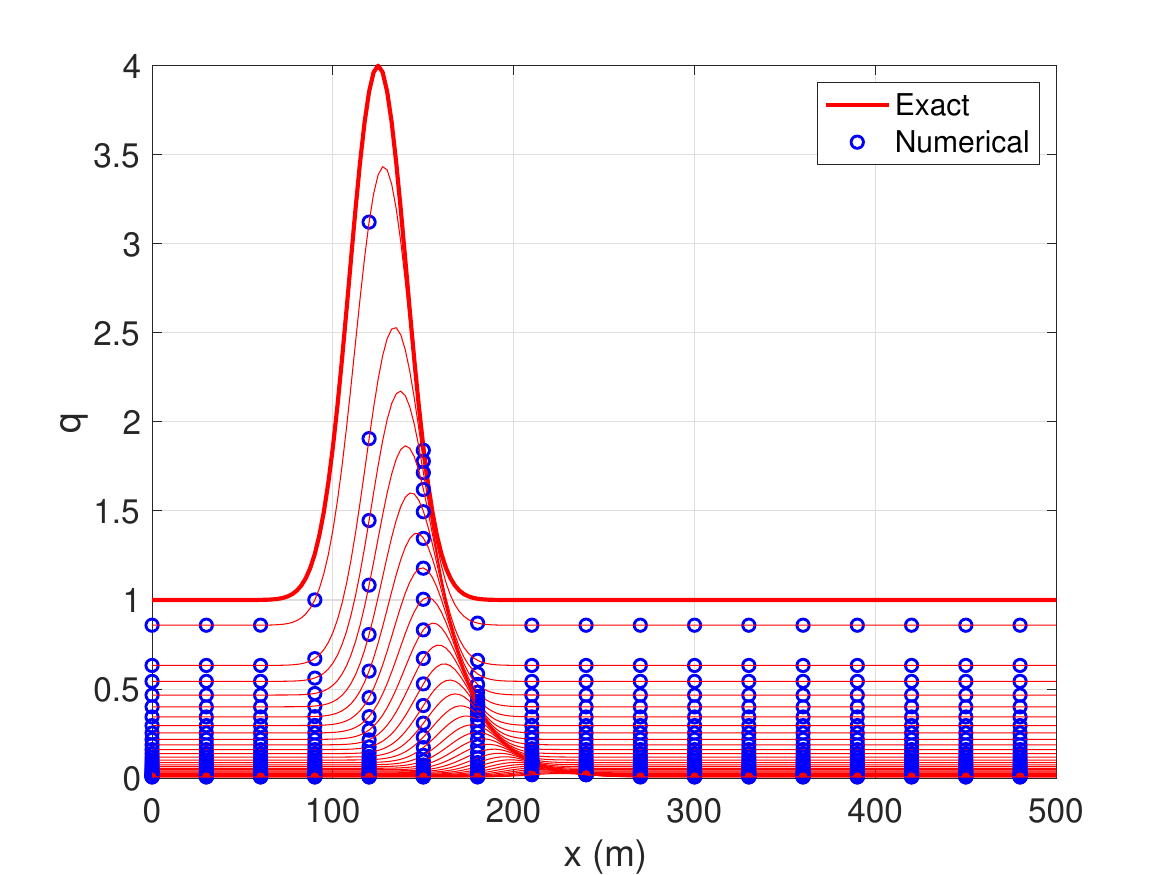}
	\caption{Linear advection-reaction problem, comparison between the exact solution (red lines) and numerical results (blue lines). Left: convergence test. Right: exact and approximate solutions over time for the mesh with $200$ elements. }
	\label{fig:linear_advectionreaction_2}
\end{figure}

Next, we focus on the linear advection-reaction model~\eqref{eq:test_equation_von_Neumann}. In this case, the exact solution is given by
\begin{equation}
	q(x, t; a; \chi) = q_{0}(x - a t)e^{-\chi t},
\end{equation}
where the initial condition is a space-varying function $q_{0} = q_{0}(x)$. We consider a smooth initial condition given by the following Gaussian bell function
\begin{equation}
	q_{0}(x) = 1 + 3\exp{\rpth{-5\rpth{\dfrac{x - L/4}{0.1L}}^{2}}},
\end{equation} 
with $L = \SI{500}{\meter}$. The final time is set to $T_{f} = \SI{100}{\second}$. In all simulations, we select the time step based on the CFL condition $\nu = \frac{a\Delta t}{\Delta x} \approx 1.22$, which is the stability limit found in the previous section for the scheme in Table~\ref{tab:IMEX_maximized_stability}, and we compute the space-time error in the $L^{\infty}$ norm. We focus on the case where we have $a=1$, $\chi=0.05$. To assess convergence, we exploit the availability of the exact solution. Simulations are run up to the final time on a sequence of uniform meshes with 100 to 500 cells. The scheme exhibits the expected second-order accuracy (Figure~\ref{fig:linear_advectionreaction_2}).

\begin{figure}[h!]
	\centering
	\includegraphics[scale = 0.5]{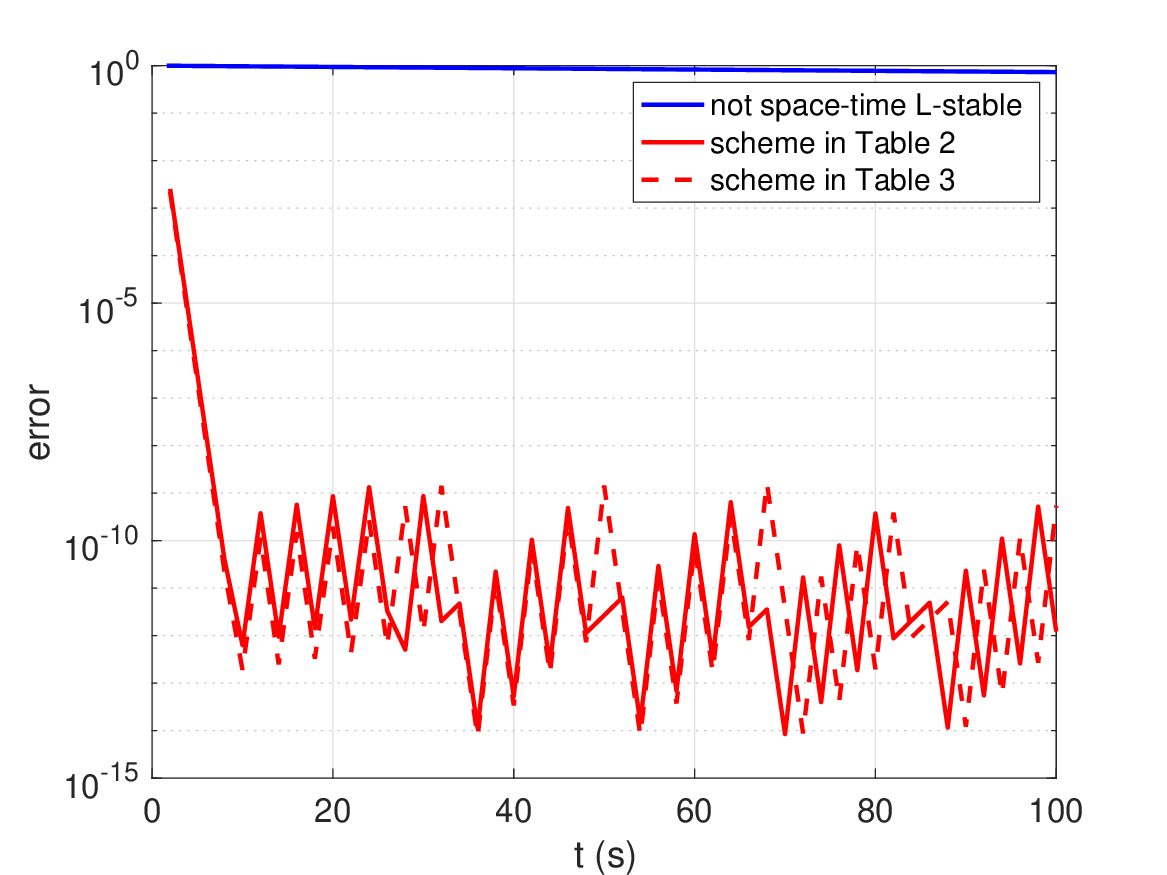}
	\caption{Linear advection-reaction problem. Error in $L^{\infty}$ norm as function of time.}
	\label{fig:linear_advectionreaction}
\end{figure}

In a second example, to test the difference with respect to the numerical scheme adopted in~\cite{gatti:2025}, motivated by the fact that the amplification factor given in Equation (A.7) of~\cite{gatti:2025} grows for higher spatial frequencies, we consider the following initial condition in the discrete one-dimensional domain $\cpth{x_{i}}_{0}^{N}$, where $N$ is the size of the discrete domain
%.
\begin{equation}
	q_{0}(x_{i}) = (-1)^{i}\exp\rpth{-\rpth{\dfrac{x_{i} - L/2}{0.05L\sqrt{2}}}^{2}}.
\end{equation}
For this test, we use the same final time as before and a mesh consisting of $300$ intervals. In all simulations, the time step is selected according to the CFL condition. We take a Courant number $\nu$ of $1.22$ for the proposed schemes and $0.9$ for the IMEX scheme used in~\cite{gatti:2025}. The advection speed is set to $a = 1$ and the reaction coefficient to $\chi = 10^{3}$, since we are interested in the asymptotic behavior. The goal is to illustrate the difference between a space-time \texttt{L}-stable method and a non-space-time \texttt{L}-stable method. In Figure~\ref{fig:linear_advectionreaction}, we report the error in the $L^{\infty}$ norm as a function of time for both the scheme proposed in~\cite{gatti:2025} and the new space-time \texttt{L}-stable schemes developed in this work. As can be observed, while the former scheme fails to adapt to the dynamics induced by the reaction term, the proposed schemes successfully capture the dynamics, yielding an approximation error that remains close to numerical zero as time advances.

\subsection{Traveling vortex}
\label{ssec:vortex}

Following the work~\cite{arpaia:2026}, we provide a numerical experiment aimed at testing the convergence of the scheme developed in a case of interest in geophysical flows. The benchmark is inspired by the compact support traveling vortex for the shallow water equations described in detail in~\cite{ricchiuto:2021}. The friction term and the contribution of the lava discharge are not included, i.e. $\lambda = Q = 0$ in~\eqref{eq:lava_model_Z}. Moreover, a flat topography is employed, i.e. $Z = 0$. An analytical solution is available, which is used therefore to assess the convergence properties of our novel numerical method and to validate our implementation. The initial conditions read as follows
\begin{subequations}
	\begin{align}
		h(r) &= h_{0} - \begin{cases}
			\frac{4}{g}\rpth{\frac{2^{p}\Gamma r_{0}}{\pi}}^{2}\rpth{H_p(\pi/2) - H_p(\rho/2)} & r\le r_{0}, \\
			0 & \text{otherwise},
		\end{cases}
		\quad \rho = \frac{\pi r}{r_{0}}, \\
		u &= u_{\infty} - \rpth{y - y_{0}}\omega(\rho), \\
		v &= \rpth{x - x_{0}}\omega(\rho).
	\end{align}
\end{subequations}
Here

$$\omega(\rho) = 
\begin{cases} 
	2^{p}\Gamma\cos^{2p}\rpth{\rho/2} & r\le r_{0}, \\
	0 & \text{otherwise},
\end{cases}
$$
where $\Gamma$ is a free parameter that controls the strength of the vortex and imposes a minimal depth $h(0) = h_{\mathrm{min}}$ in the center of the vortex $\rpth{x_{0}, y_{0}}$, while $p$ is an integer parameter that controls the regularity of the solution. More specifically, the vortex is $C^{2p}(\Omega)$. Moreover, $r = \sqrt{\rpth{x - x_{0}}^{2} + \rpth{y - y_{0}}^{2}}$ denotes the distance from the center of the vortex, while $u_{\infty}$ is the background velocity oriented along the $x$ axis at which the vortex is transported. The definition of $H(x)$ is the following
\begin{align*}
	H_{p}(x) &= \int x\cos^{4p}(x)\mathrm{d}x 
	= \frac{x\cos^{4p-3}(x)\sin(x)}{4p}\spth{\cos^{2}(x) + \frac{4p - 1}{4p - 2}} \\&+ \cos^{4p-2}(x)\spth{\frac{4p-1}{4p\rpth{4p-2}^{2}} + \frac{\cos^{2}(x)}{\rpth{4p}^{2}}} 
	+ \frac{4p-1}{4p}\frac{4p-3}{4p-2}\int x\cos^{4p - 4}(x)\mathrm{d}x.
\end{align*}
We consider $p = 1$, so that we obtain the following expressions 
$$H_{1}(x) = \frac{x\cos(x)\sin(x)}{4}\spth{\cos^{2}(x) + \frac{3}{2}} + \cos^{2}(x)\spth{\frac{3 + \cos^{2}(x)}{16}} + \frac{3}{16}x^{2}$$
and
$$\Gamma = \frac{\pi}{2 r_{0}}\sqrt{\frac{g(h_{0} - h_{\mathrm{min}})}{H_{1}(\pi/2) - H_{1}(0)}}.$$
We recall that in the formula above $g = \SI{9.81}{\meter\per\second\squared}$ denotes the acceleration of gravity. The parameter values defining the vortex and the background flow are the following
$$h_{0} = \SI{1}{\meter}, \quad h_{\mathrm{min}} = \SI{0.9}{\meter}, \quad u_{\infty} = \SI{6}{\meter\per\second}, \quad \rpth{x_{0}, y_{0}} = \rpth{\SI{0.5}{\meter},\SI{0.5}{\meter}}, \quad r_{0} = \SI{0.25}{\meter}.$$
The final time is $T_{f} = \SI[parse-numbers=false]{1/6}{\second}$, while the computational domain is given by $(0,2)\SI{}{\meter} \times (0,1)\SI{}{\meter}$. A constant temperature $T = \SI{300}{\kelvin}$ is employed. The exact solution is a propagation of the initial condition at the background velocity.

\begin{figure}[h!]
	\centering
	\includegraphics[scale = 0.5]{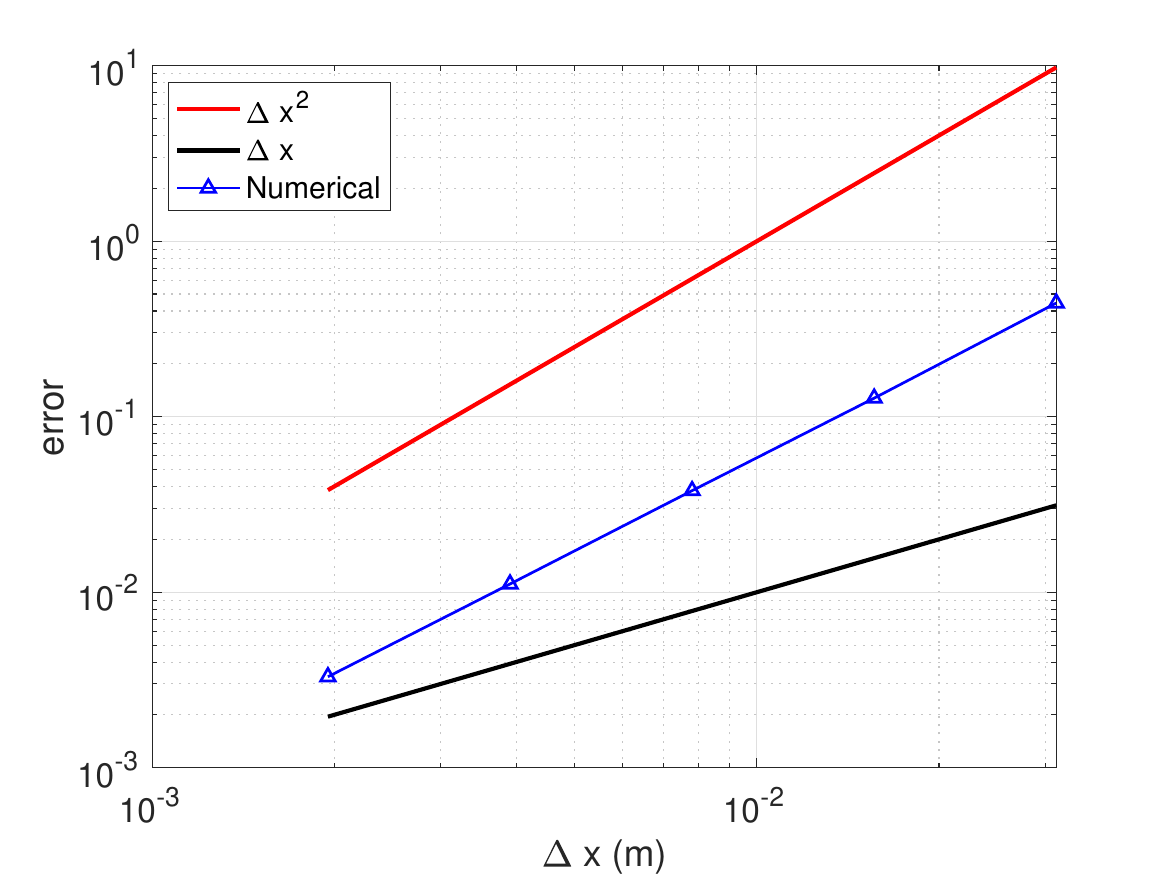}
	\caption{Traveling vortex. Convergence plot. The black line shows the first order convergence, while the red line displays the second order one. The blue line shows the numerical order of convergence.}
	\label{fig:linear_travel}
\end{figure}

We propose a convergence analysis to verify whether the developed method is able to achieve second-order convergence when nonlinear convective terms are present. For this test, we consider the Courant number fixed at $1.1$. Since both methods developed in Table~\ref{tab:IMEX_c_eq_ctilde} and Table~\ref{tab:IMEX_maximized_stability} behave in a similar way, in the following we present the results obtained only with the method of Table~\ref{tab:IMEX_maximized_stability}. To perform the convergence test, we consider a set of structured meshes with a level of refinement ranging from $6$ to $10$. In Figure~\ref{fig:linear_travel}, we report the convergence plot obtained considering only the numerical approximation along the centerline, i.e., $y = \SI{0.5}{\meter}$, since the expected dynamics is purely one-dimensional. The errors are computed in the $L^{2}$ norm. As can be seen, the method is close to second-order convergence; in particular, we obtain an estimated convergence rate of approximately $1.77$.

\begin{figure}[h!]
	\centering
	\includegraphics[scale = 0.2]{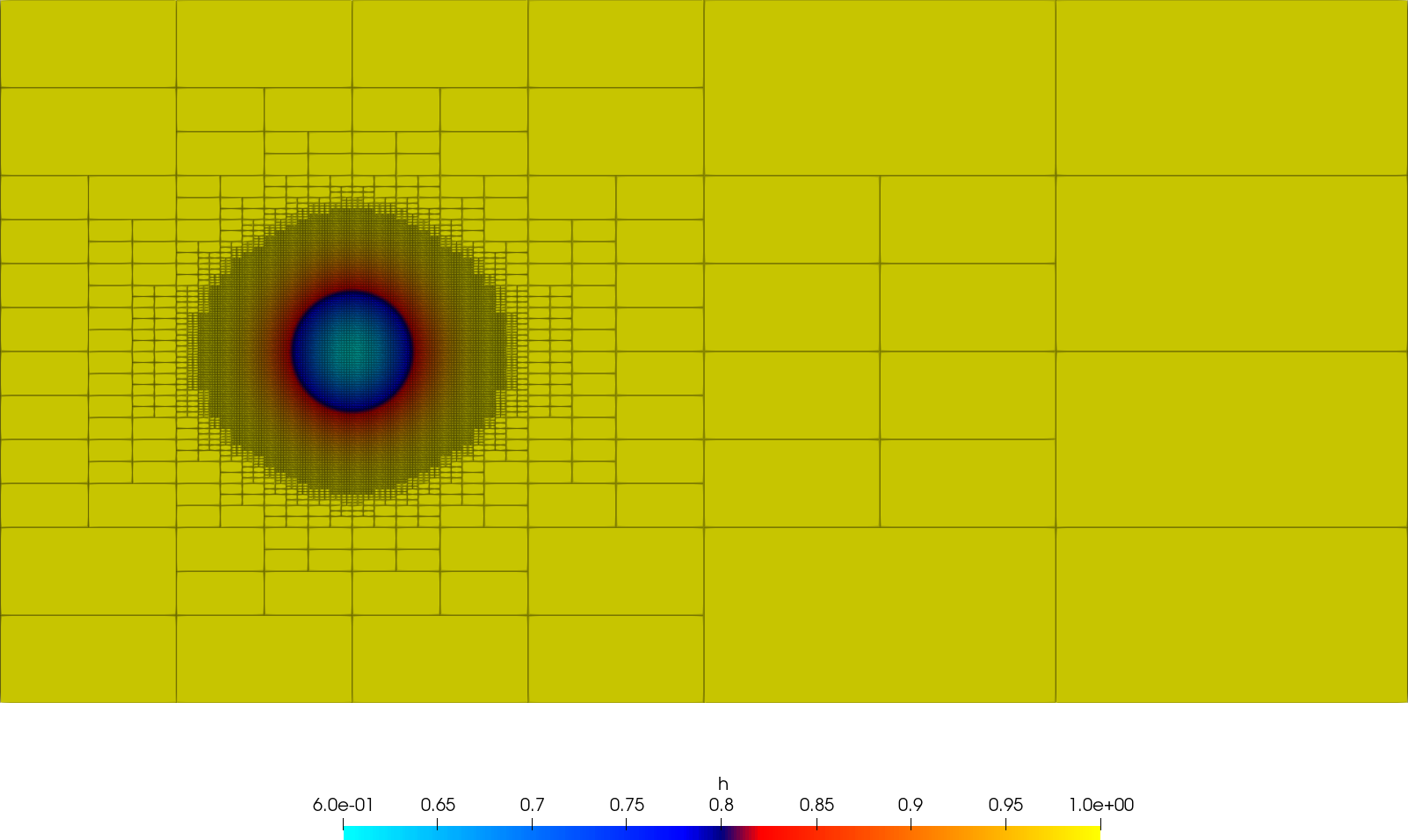}\vfill
	\includegraphics[scale = 0.5]{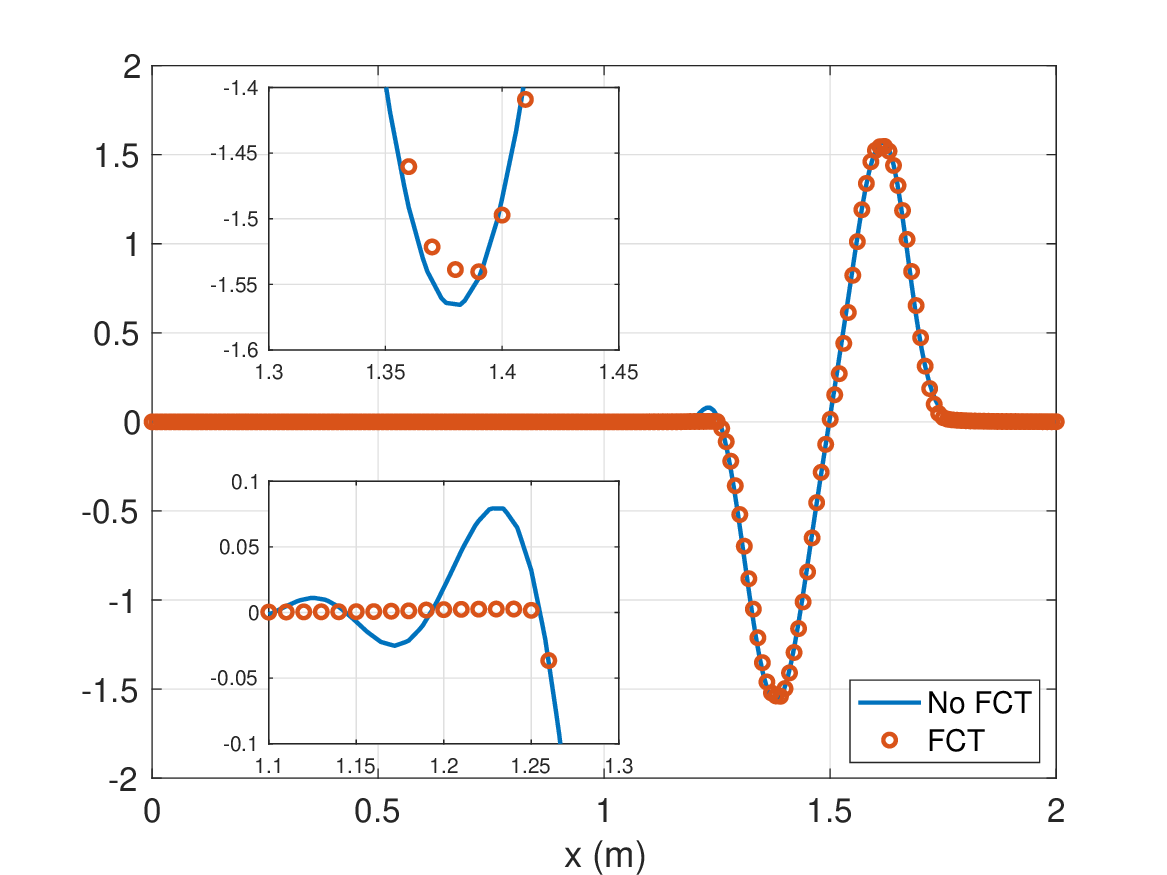}
	\caption{Traveling vortex. Top: Sketch of the isolines of the initial material height together with the initial mesh in case the maximum level of refinement present is $9$. Bottom: effect of the FCT procedure (red dots) on the proposed scheme.}
	\label{fig:linear_travel_im}
\end{figure}

Despite the smoothness of the problem, we notice that the adopted method shows spurious oscillations in the approximation of the mass flux. Numerical oscillations of the same kind also appear when using a classical Lax--Wendroff procedure in line with the one described in~\cite{gatti:2023b}. We now consider an adaptive mesh quadtree framework, and we refer to~\cite{gatti:2023b} for all the details concerning the mesh adaptation procedure, the estimator requires a tolerance that we set equal to $\SI{0.01}{\meter}$. We show in Figure~\ref{fig:linear_travel_im} the isolines of the initial material height together with the quadtree mesh superimposed on it. The maximum level of mesh refinement that we set is equal to $9$. Moreover, similar to the approach described in~\cite{gatti:2023b}, we apply the Flux Corrected Transport (FCT) procedure as an artificial diffusion strategy to limit the reconstruction of the higher-order fluxes that result from the evaluation of $\mathbf{q}^{(n,3)}$ in~\eqref{eq:eq_stage_3_up}. We provide a comparison between the solution obtained when the FCT procedure is applied and when it is not (Figure~\ref{fig:linear_travel_im}).
In particular, we focus on the one-dimensional result obtained along the line $y = \SI{0.5}{\meter}$ and show only the mass flux along the dimension $y$ $hu_{y}$, as we experience the rise of oscillations in that solution. As one can notice, for the present case the FCT procedure is able to dissipate spurious oscillations when applied in conjunction with the developed scheme. One can also notice that the adopted FCT procedure is unable to maintain local maxima, as one can see from one of the two zoomed windows in Figure~\ref{fig:linear_travel_im}.

We note that, for the newly developed schemes proposed in this work, the conditions under which a flux-limiting procedure such as FCT is required, when applied to the reconstruction of $\mathbf{q}^{(n,3)}$, remain to be fully clarified. In the present test case, a clear benefit is observed, as the procedure effectively damps spurious oscillations to some extent. However, in general, in the presence of reaction contributions, the use of an FCT procedure may not be strictly necessary and depends on the specific problem, since the reaction term itself can provide an intrinsic smoothing effect, thereby reducing the need for additional artificial diffusion.For this reason, we do not employ any FCT procedure in the subsequent numerical experiments, in particular those of Section~\ref{ssec:pouring_lava}, as a stiff reaction term is present. The preservation of physical constraints, which is closely related to the control of spurious oscillations, is deferred to future investigations.

\subsection{Well-balancing tests}
\label{ssec:wb_tests}

\begin{figure}[h!]
	\centering
	\includegraphics[scale = 0.15]{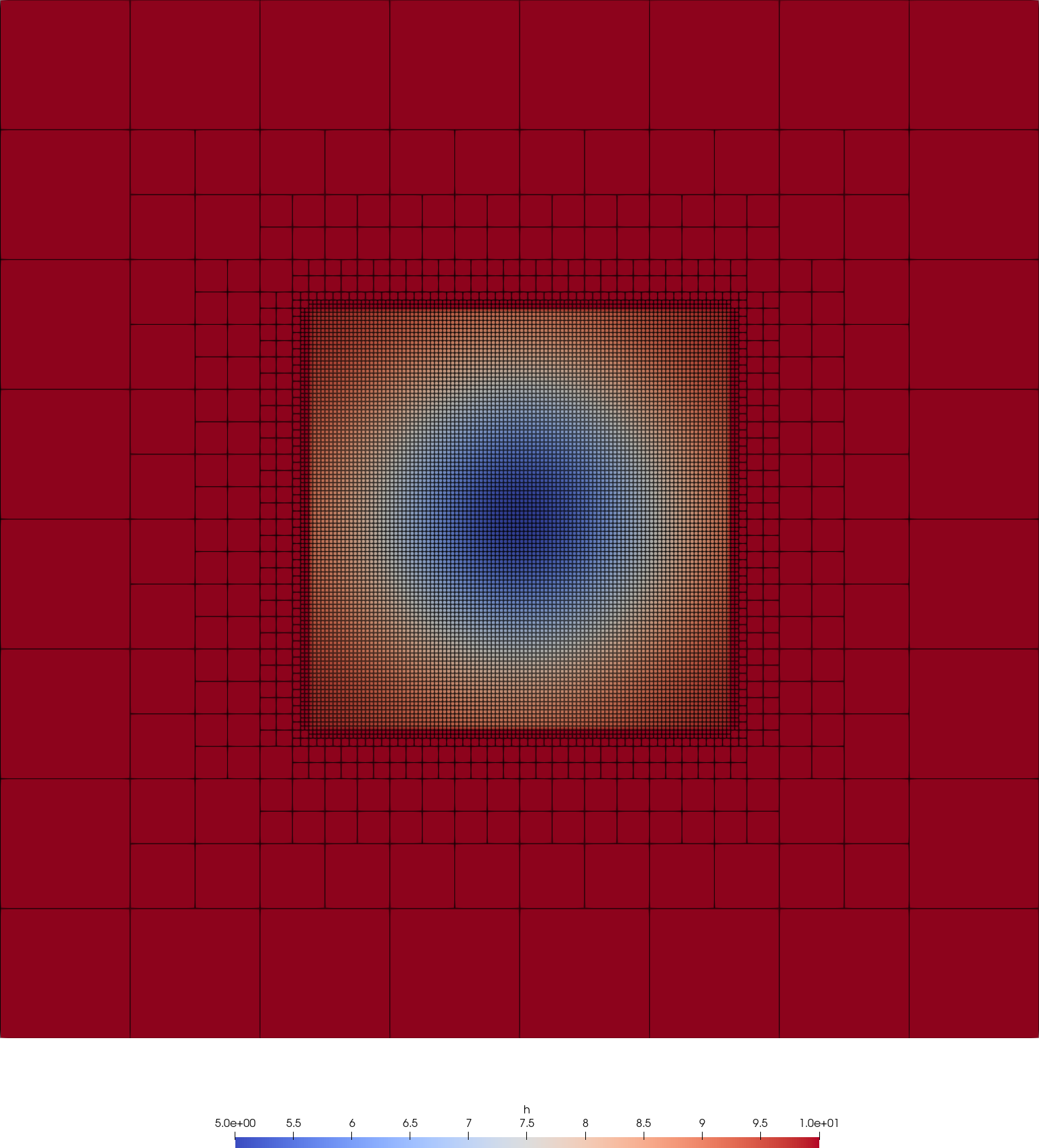}
	\caption{Well-balancing tests. Sketch of the initial material height together with the employed quadtree mesh superimposed on it.}
	\label{fig:wb}
\end{figure}

Next, we verify the well-balancing property of the numerical method. We analyze an example that blends both smooth and non-smooth topographies inside it. We consider $L_{x} = L_{y} = \SI{10}{\meter}$ and $T_{f} = \SI{1}{\second}$. Wall boundary conditions are prescribed. For what concerns the initial conditions, we assume the lava-at-rest condition, with a total free surface $\zeta = \SI{10}{\meter}$ and $hT = \SI[parse-numbers=false]{10^{3}}{\meter\kelvin}$. The topography is defined by
\begin{equation}
	Z(\bm{x}) = 
	\begin{cases}
		5\exp\spth{-\frac{2}{5}\rpth{x-5}^{2} - \frac{2}{5}\rpth{y-5}^{2}} & \text{if $3\le x \le 7$ and $3 \le y \le7$}, \\
		0 & \text{otherwise}.
	\end{cases}
\end{equation}
We consider an adaptive mesh built with a tolerance $\SI{0.001}{\meter}$ and a minimum mesh size equal to $\SI{0.05}{\meter}$ (see Figure~\ref{fig:wb}), the time step is chosen by considering a Courant number equal to $1.1$. We report in Table~\ref{tab:linf_error_WB} the errors in the $L^{\infty}$ norm at the final time. We note that we perform roughly $280$ time steps. As the reader can see, the method is able to preserve the lava-at-rest state since the errors are close to the round-off unit.
\begin{table}[h!]
	\centering
	\def\arraystretch{1.3}\tabcolsep=8pt
	\begin{tabular}{c c c c c}
		\hline
		& $h$ & $hu_{x}$ & $hu_{y}$ & $hT$ \\ 
		\hline
		$Z$ & $7.10\times 10^{-15}$ & $7.74\times 10^{-14}$ & $7.44\times 10^{-14}$ & $9.99\times 10^{-16}$ \\
		\hline
	\end{tabular}
	\caption{Well-balancing tests. $L^{\infty}$ norm of the error for the lava-at-rest solution.}
	\label{tab:linf_error_WB}
\end{table}

\subsection{Pouring of lava from a vent over a flat topography}
\label{ssec:pouring_lava}

In this final test, we simulate the pouring of lava from a vent over a flat topography under considerably stiff source conditions. The goal is to demonstrate that the novel second-order (pseudo-)staggered Galerkin method developed in this work can successfully handle cases in which the IMEX schemes of~\cite{gatti:2025} fail. The domain is $\Omega = (0,200) \times(0,200)\,\si{\meter\squared}$, with the vent located in the center, $x_{v} = y_{v} = \SI{100}{\meter}$ and characterized by $\sigma = \SI[parse-numbers=false]{10^{-1}}{\meter\squared}$ in~\eqref{eq:Gaussian_lava}. The lava exits the vent at a temperature of $T_{e} = \SI{2000}{\kelvin}$ with a constant discharge of $Q = \SI{200}{\cubic\meter\per\second}$. The coefficient $b$ in the exponential stiff source term~\eqref{eq:lambda_vent} is set to $b = \SI[parse-numbers=false]{10^{-2}}{\per\kelvin}$. The simulations are performed with a space adaptation procedure, using a tolerance of $\SI[parse-numbers=false]{10^{-2}}{\meter}$, updated every $\SI{0.1}{\second}$, and a minimum mesh size of $\SI{0.1}{\meter}$. A saturation threshold of $10^{8}$ is imposed on the stiff source term in the momentum equation. As at the beginning of the simulation there is no lava in the domain, we just set the initial time step equal to $\SI[parse-numbers=false]{10^{-4}}{\second}$.

\begin{figure}[h!]
	\centering
	\includegraphics[scale = 0.5]{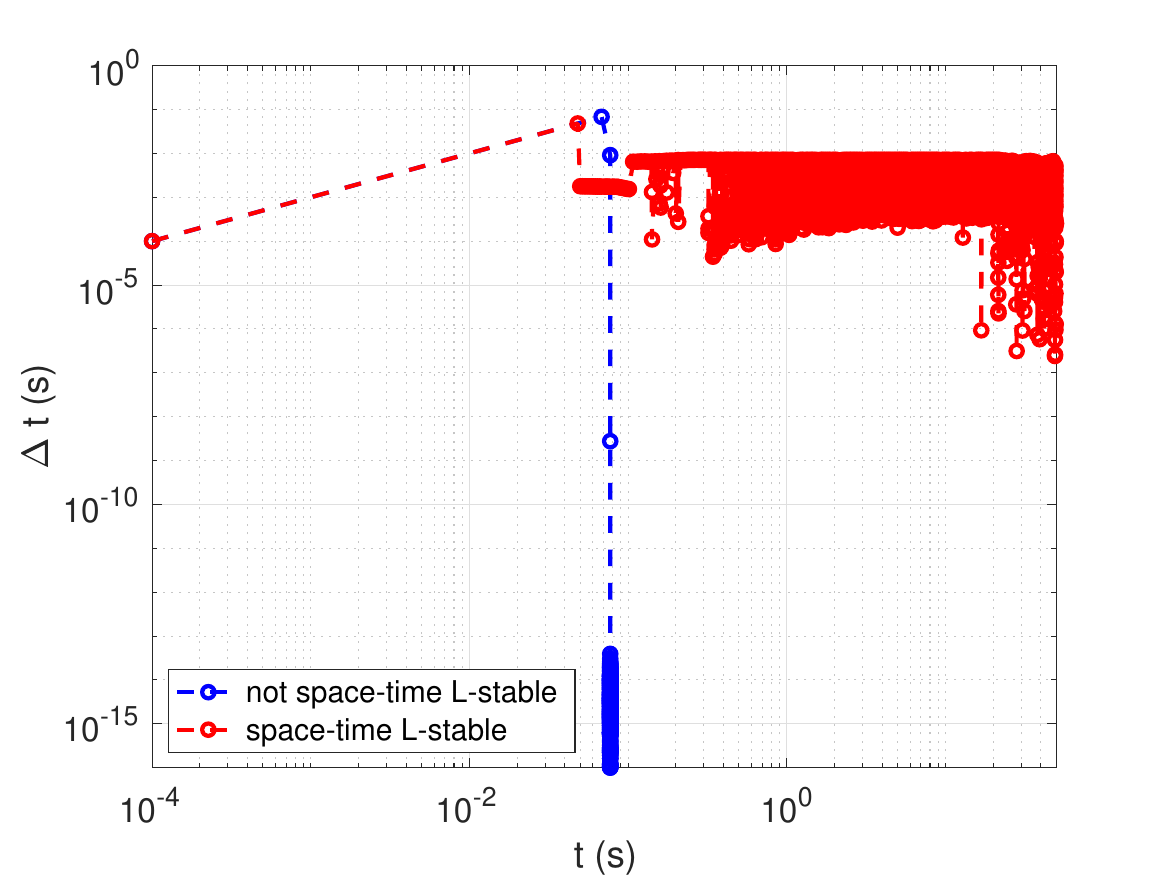}
	\caption{Pouring of lava from a vent over a flat topography. Time step size over simulation time for the method developed in~\cite{gatti:2025} (blue dots) and the space-time \texttt{L}-stable discretization in Table~\ref{tab:IMEX_maximized_stability} (red dots).}
	\label{fig:lava_pouring_time_step}
\end{figure}

The comparison between the newly developed schemes and that reported in~\cite{gatti:2025} is striking. We use a Courant number of $0.9$ for the method in~\cite{gatti:2025} and $1.1$ for the space-time \texttt{L}-stable schemes introduced in this work. Figure~\ref{fig:lava_pouring_time_step} shows the evolution of the time step size for both methods. The method used in~\cite{gatti:2025} cannot handle the stiffness of the problem, and the time step rapidly approaches the machine epsilon. By contrast, the (pseudo)staggered Galerkin method with IMEX-RK coefficients in Table~\ref{tab:IMEX_c_eq_ctilde} and Table~\ref{tab:IMEX_maximized_stability} exhibits a stable solution, thanks to their intrinsic space-time \texttt{L}-stability, which effectively filters out transient dynamics induced by the stiff source term. Without this property, the time step would eventually collapse to zero. We just report the time step size over time for only one of the two proposed methods, as both exhibit similar performance.

\begin{figure}[h!]
	\centering
	\includegraphics[scale = 0.16]{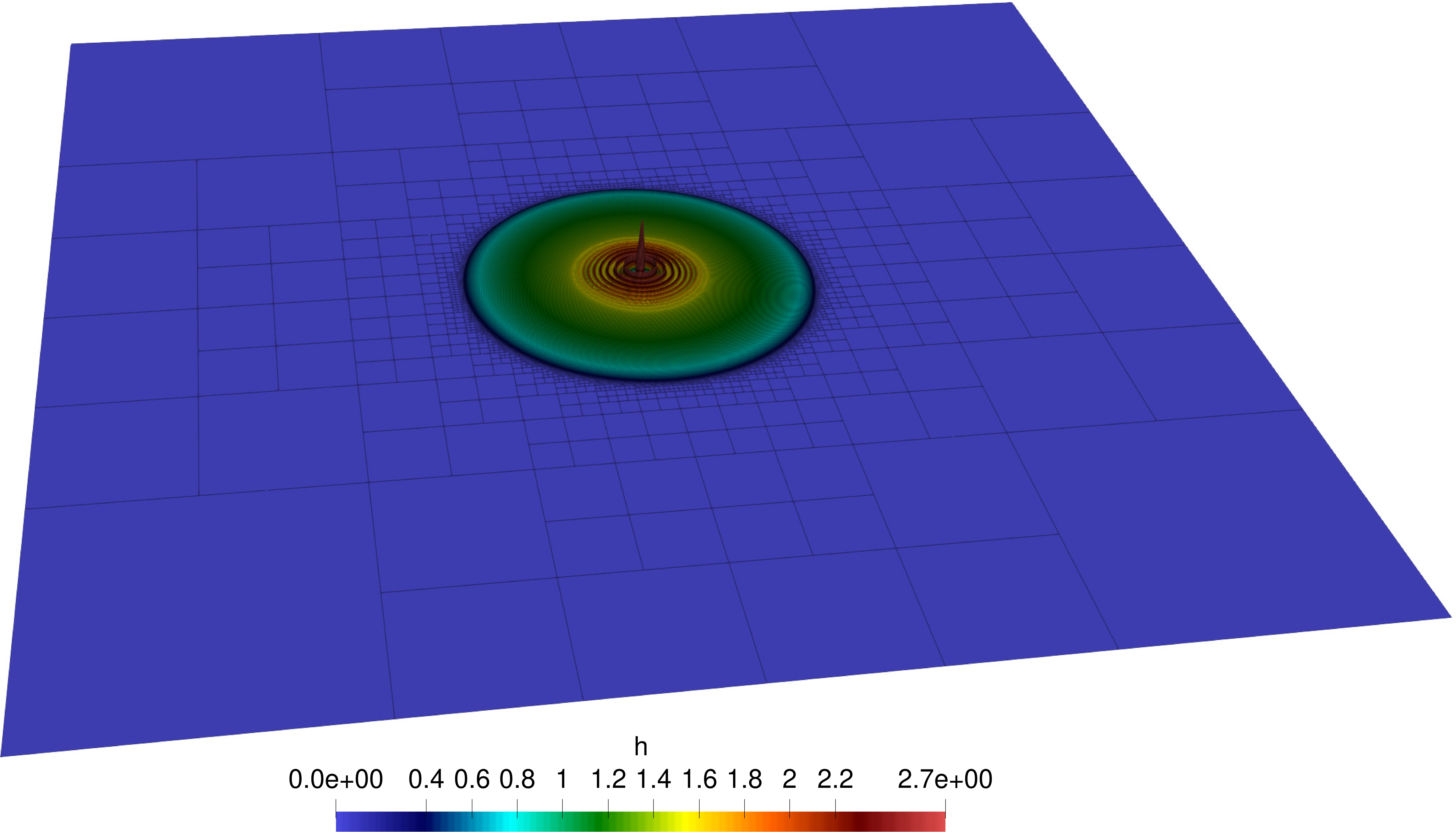}
	\includegraphics[scale = 0.5]{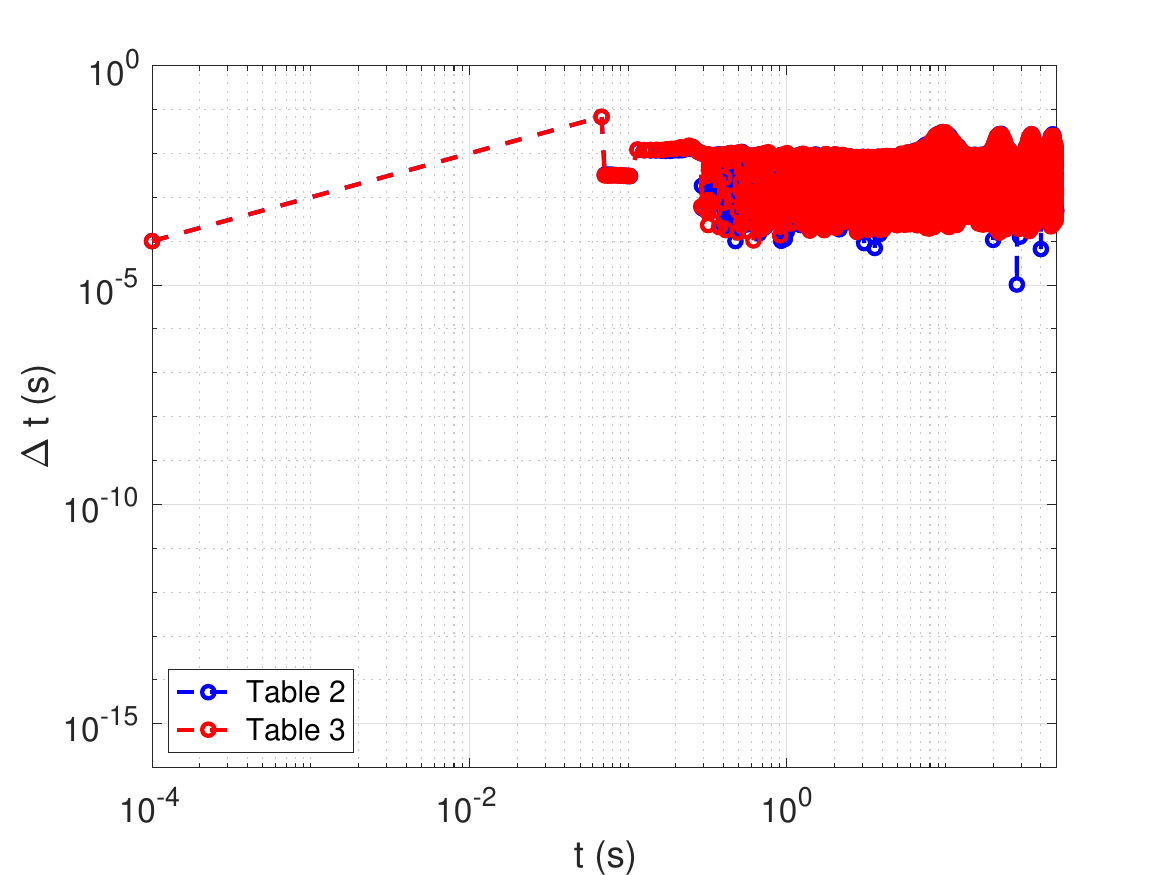}
	\caption{Pouring of lava from a vent over a flat topography. Top: Sketch of the approximation of the material height at $\SI{90}{\second}$ of simulation time together with the corresponding quadtree mesh. Bottom: evolution of the time step for the IMEX-RK scheme in Table~\ref{tab:IMEX_c_eq_ctilde} (blue dots) and in Table~\ref{tab:IMEX_maximized_stability} (red dots).}
	\label{fig:lava_pouring_time_step_2}
\end{figure}

For the sake of completeness, we now consider a time-dependent lava discharge modeled as a nonlinear, modulated chaotic map designed to produce aperiodic oscillations. More specifically, we set
\begin{equation}
	Q(t) = \max\rpth{0,\, \tilde{Q}\rpth{\frac{1}{2} + \frac{1}{2} \sin(4t)}\cos(12t)},
\end{equation}
where $\tilde{Q}$ is obtained by iterating a logistic map three times,
\begin{equation}
	\tilde{Q} = f\rpth{f\rpth{f\rpth{Q_{0}}}}, \qquad f(x) = r \, x \, (1 - x), \qquad r = 3.7,
\end{equation}
with $Q_{0} = \SI{200}{\cubic\meter\per\second}$ denoting the initial value. This construction combines the chaotic dynamics of the logistic map with deterministic oscillatory modulation (th\-rough $\sin(4t)$ and $\cos(12t)$), producing a complex aperiodic forcing signal with high- and low-frequency components. The non-negativity condition enforces that we are dealing with a source vent and not with a sink point.

For this time-dependent forcing term as well, both sets of coefficients listed in Table~\ref{tab:IMEX_c_eq_ctilde} and Table~\ref{tab:IMEX_maximized_stability} successfully complete the simulation. However, the coefficients in Table~\ref{tab:IMEX_c_eq_ctilde} lead to smaller time steps at certain phases of the simulation (see Figure~\ref{fig:lava_pouring_time_step_2}, bottom). This behavior can be explained by noting that the scheme in Table~\ref{tab:IMEX_maximized_stability} is derived by maximizing the stability region in the advection-dominated regime. However, despite the generally larger stability region of the latter scheme under advection-dominated conditions, the scheme in Table~\ref{tab:IMEX_c_eq_ctilde} requires fewer total time steps than that in Table~\ref{tab:IMEX_maximized_stability}: specifically, $25{,}418$ time steps are required by the scheme in Table~\ref{tab:IMEX_c_eq_ctilde}, compared to the $34{,}500$ time steps required by the method in Table~\ref{tab:IMEX_maximized_stability}. 
We remark that the smaller time steps observed for the scheme in Table~\ref{tab:IMEX_c_eq_ctilde} are confined to specific phases of the simulation, presumably when the oscillatory forcing $Q(t)$ most strongly excites high-frequency components of the solution. Outside these intervals, this scheme allows for systematically larger time steps than the scheme in Table~\ref{tab:IMEX_maximized_stability}, which instead enforces a more uniformly conservative time step selection throughout the simulation. This trade-off explains why, despite its locally reduced stability region in the advection-dominated regime, the scheme in Table~\ref{tab:IMEX_c_eq_ctilde} ultimately requires fewer total time steps ($25{,}418$) to complete the simulation than the scheme in Table~\ref{tab:IMEX_maximized_stability} ($34{,}500$), which in turn trades computational efficiency for increased robustness against oscillatory forcing.

We report in Figure~\ref{fig:lava_pouring_time_step_2} (top) the approximation of the material height after $\SI{90}{\second}$ of the simulation time. We note that the solution is qualitatively similar to the one displayed in~\cite{gatti:2025}, but it presents many more spatial oscillations due to the heavy stiffness and the aperiodic oscillations of the lava discharge. Moreover, the spatial symmetry of the space adaptation procedure is preserved.

%%%%%%%%%% Conclusions %%%%%%%%%%%%%%%%%%%%%
\section{Conclusions and further developments}
\label{sec:conclu}

We have proposed a novel sec\-ond-order (pseudo-)staggered Galerkin method which is particularly well suited for the solution of a modified set of shallow water equations that model the dynamics of lava flows. The time discretization is based on an Implicit-Explicit (IMEX) Runge--Kutta (RK) method, whose coefficients have been determined following the outcome of a von Neumann stability analysis. In particular, the IMEX (pseudo)staggered Galerkin method satisfies a space-time \texttt{L}-stability condition, which is useful for avoiding oscillations in the mass flux. Moreover, a Lax--Wendroff procedure has been performed to further enhance the stability properties of the numerical method. The accuracy and well-balancing properties for the lava flow model of this novel numerical method have been assessed in a number of relevant test cases.

In future work, the authors are planning to perform a quantitative comparison between this novel (pseudo-)staggered Galerkin method and a discretization scheme based on the Discontinuous Galerkin (DG) method. 
The DG method has been widely used to develop numerical methods for SWE~\cite{arpaia:2026, dumbser:2013, tumolo:2013}, for the lava flow model~\cite{conroy:2021}, and, more generally, for environmental flows; see, among others~\cite{orlando:2022, orlando:2023, orlando:2024, tumolo:2013}. One goal is to mimic the von Neumann stability analysis performed in this work also for DG discretization. However, because of the locality of the DG method, the space-time \texttt{L}-stability condition introduced in this work reduces to the classical \texttt{L}-stability property for the discretization of ordinary differential equations. Hence, a larger number of free parameters remains for the IMEX scheme, whose value can be determined, e.g., according to an absolute monotonicity analysis or requiring a high-order of dissipation or dispersion. The derivation of some second-order IMEX schemes in combination with the DG method and their comparison with the results of this article will be the main topic of future work.

A related question concerns the extension of the present framework to higher-order IMEX methods. In that case, the number of free parameters also increases significantly, and the information provided by the von Neumann stability analysis alone may not be sufficient to fully determine the implicit companion method, unlike the three-stage second-order schemes considered in this work. It would therefore be interesting to investigate whether additional design principles, analogous to those required in the DG setting, can be employed to select the remaining degrees of freedom. Nevertheless, the analysis carried out in this work identifies a set of stability requirements arising from the interplay between temporal and spatial discretizations. To the best of the authors' knowledge, exploiting such requirements as design criteria for IMEX schemes is novel and constitutes one of the main contributions of the present paper.

\section*{Acknowledgements}
G.O. is part of the INdAM-GNCS National Research Group. Generative AI tools have been use exclusively to improve language clarity, without contributing to the scientific content of the work. The authors assume responsibility for all content.

\subsection*{Conflict of interest}
The authors declare that they have no potential conflict of interest.

\subsection*{Code availability}

\noindent
{The numerical method described in this article has been implemented in C++ and can be accessed through the following dedicated repository on GitHub: \url{https://github.com/federicg/lava-flow}. The \texttt{scripts} folder of the repository contains two Mathematica notebooks, \texttt{LW.nb} for the Lax--Wendroff procedure and \texttt{von\_Neumann.nb} for the von Neumann stability analysis, as well as a MATLAB script providing a one-dimensional implementation of the numerical method presented in this work.}

%%%%%%%%%%%%% Appendices %%%%%%%%%%%%%%%%%%%%
\appendix

\section{Computation of the non-conservative term in the bilinear space}
\label{app:detail_wb}

\begin{figure}[h!]
    \centering
    \begin{tikzpicture}[scale=2]
        \draw[green, very thick] (0,0) rectangle (3,3);
        \fill[red] (0,0) circle[radius=3pt];
        \fill[red] (3,3) circle[radius=3pt];
        \fill[red] (0,3) circle[radius=3pt];
        \fill[red] (3,0) circle[radius=3pt];
        \path (1.5,1.5) pic[black,very thick] {cross=3pt};
        \node[red, below left=1pt, outer sep=2pt] at (0,0) {$(i,j)$};
        \node[red, below right=1pt, outer sep=2pt] at (3,0) {$(i+1,j)$};
        \node[red, above left=1pt, outer sep=2pt] at (0,3) {$(i,j+1)$};
        \node[red, above right=1pt, outer sep=2pt] at (3,3) {$(i+1,j+1)$};
        \node[above=1pt, outer sep=2pt] at (1.5,1.5) {$(i+\tfrac{1}{2},j+\tfrac{1}{2})$};
    \end{tikzpicture}
    \caption{Quadrilateral element centered at $(i+\frac{1}{2},j+\frac{1}{2})$. The degrees of freedom associated with space $\mathbb{Q}_{1}$ are red-highlighted; the black cross (i.e., the barycenter of the quad) corresponds to the degree of freedom of the space $\mathbb{Q}_{0}$.}
    \label{fig:appendix1}
\end{figure}
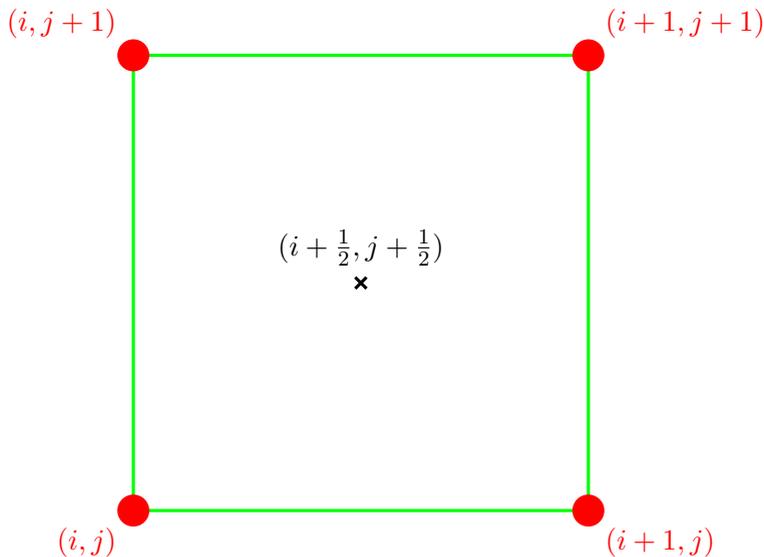

Here, we clarify the computation of the non-conservative contributions in the bilinear space $\tilde{\bm{Q}}_{1}$. These contributions appear both in the first step of the method and in the update step, whereas for the computation of the solution stage $\mathbf{q}^{(n,3)}$ we rely on the PC method. It is important to follow the prescribed procedure carefully in order not to deteriorate the well-balancing property, which is guaranteed when this procedure is combined with the PC method (see~\cite{gatti:2024a} for details on the implementation of the PC method). \\
We focus on a quadrilateral mesh element of size $\Delta x \times \Delta y$, as shown in Figure \ref{fig:appendix1}. Since the same procedure applies to both the intermediate step $\mathbf{q}^{(n,2)}$ and the updated solution $\mathbf{q}^{n+1}$, we show only the explicit computation for the former case. For the updated solution, it is only necessary to note that the solution is defined on the mesh nodes rather than at the mesh centroid. In this case, the solution update is obtained from the sum of the contributions of all mesh quadrants sharing the same node.

During the computation of the intermediate solution $\mathbf{q}^{(n,2)}$, the integration of the pressure inertial contribution to the mass flux $hu_{x}$ produces
\begin{align}
    \rpth{\dfrac{1}{2}g\partial_{x}h^{2}, \phi^{(0)}_{i+\frac{1}{2}, j+\frac{1}{2}}} &= \rpth{\dfrac{1}{2}\dfrac{\tfrac{1}{2}g h^{2}_{i+1,j} - \tfrac{1}{2}g h^{2}_{i,j}}{\Delta x} + \dfrac{1}{2}\dfrac{\tfrac{1}{2}g h^{2}_{i+1,j+1} - \tfrac{1}{2}g h^{2}_{i,j+1}}{\Delta x}}\Delta x\Delta y \nonumber \\[2mm]
    &= \dfrac{1}{2}\left(g\dfrac{h_{i,j} + h_{i+1,j}}{2}\dfrac{h_{i+1,j} - h_{i,j}}{\Delta x} \right. \nonumber \\ 
    &\qquad\left. + g\dfrac{h_{i,j+1} + h_{i+1,j+1}}{2}\dfrac{h_{i+1,j+1} - h_{i,j+1}}{\Delta x}\right)\Delta x\Delta y.
\end{align}
It is then clear that, to obtain a well-balanced stage and solution update, the non-conservative contributions must be integrated consistently.
From the last equation, it is evident that well-balancing is achieved when using a slope discretization such as

\begin{align}
    \rpth{gh\partial_{x}Z, \phi^{(0)}_{i + \frac{1}{2}, j + \frac{1}{2}}} 
    &= \dfrac{1}{2}\left(g\dfrac{h_{i,j} + h_{i+1,j}}{2}\dfrac{Z_{i+1,j} - Z_{i,j}}{\Delta x}\right. \nonumber \\ 
    &\qquad\left. + g\dfrac{h_{i,j+1} + h_{i+1,j+1}}{2}\dfrac{Z_{i+1,j+1} - Z_{i,j+1}}{\Delta x} \right)\Delta x \Delta y .
\end{align}
From this equation, we see that well-balancing is obtained by taking the arithmetic mean of the non-conservative contributions along the two quadrilateral edges parallel to the $x$ direction. A similar procedure applies along the $y$ direction for the force balancing of the mass flux $hu_{y}$. To reach well-balancing, in the above passages we have implicitly assumed that the currently known solution satisfies $h + Z = \text{constant}$ in the mesh nodes.

%%%%%%%%%%% Details von Neumann %%%%%%%%%%%%
\section{Details of the von Neumann stability analysis}
\label{app:staggered_von_Neumann_computations}

In this section, we report the details of the computations of the von Neumann stability analysis for \eqref{eq:test_equation_von_Neumann} which leads to \eqref{eq:Glim}. As already mentioned in Section \ref{ssec:staggered_von_Neumann_et_al}, one assumes that
$$q_{j}^{0} = \exp\rpth{ik x_{j}},$$
where we recall that $i$ represents the imaginary unit and $x_{j} = j\Delta x$. Substituting the previous expression into \eqref{eq:fully_discrete_von_Neumann_stage2}, one obtains
\begin{equation}\label{eq:amplification_von_Neumann_stage2}
    q_{j+\frac{1}{2}}^{(0,2)} = \frac{1}{1 + \gamma\Phi}\spth{\rpth{1 - \tilde{a}_{21}\Phi}\frac{1 + \exp\rpth{ik\Delta x}}{2} - a_{21}\nu\rpth{\exp\rpth{ik\Delta x} - 1}}q_{j}^{0}.
\end{equation}
Denote by $G^{(n,2)}$
\begin{equation}
    G^{(n,2)} = \frac{1}{1 + \gamma\Phi}\spth{\rpth{1 - \tilde{a}_{21}\Phi}\frac{1 + \exp\rpth{ik\Delta x}}{2} - a_{21}\nu\rpth{\exp\rpth{ik\Delta x} - 1}}
\end{equation}
and notice that
\begin{subequations}
\begin{align}
    q_{j-\frac{1}{2}}^{(0,2)} - q_{j+\frac{1}{2}}^{(0,2)} &= G^{(n,2)}\rpth{\exp\rpth{-ik\Delta x} - 1}q_{j}^{0} \label{eq:deltaq_von_Neumann_stage2} \\
    q_{j-\frac{1}{2}}^{(0,2)} + q_{j+\frac{1}{2}}^{(0,2)} &= G^{(n,2)}\rpth{\exp\rpth{-ik\Delta x} + 1}q_{j}^{0} \label{eq:sumq_von_Neumann_stage2}
\end{align}
\end{subequations}
Analogously, for what concerns the third stage of the IMEX scheme \eqref{eq:fully_discrete_von_Neumann_stage3}, one obtains
\begin{align}\label{eq:amplification_von_Neumann_stage3_tmp}
    &q_{j}^{(0,3)} = \frac{1}{1 + \gamma\Phi}\left[\rpth{1 - \tilde{a}_{31}\Phi - a_{31}\nu\frac{\exp\rpth{ik\Delta x} - \exp\rpth{-ik\Delta x}}{2}}q_{j}^{0} - \right. \nonumber \\
    &\hspace{2.8cm}\left. \tilde{a}_{32}\Phi\frac{q_{j-\frac{1}{2}}^{(n,2)} + q_{j+\frac{1}{2}}^{(n,2)}}{2} + a_{32}\nu\rpth{q_{j-\frac{1}{2}}^{(n,2)} - q_{j+\frac{1}{2}}^{(n,2)}}\right],
\end{align}
so that substituting \eqref{eq:deltaq_von_Neumann_stage2} and \eqref{eq:sumq_von_Neumann_stage2} into \eqref{eq:amplification_von_Neumann_stage3_tmp}, one obtains
\begin{align}\label{eq:amplification_von_Neumann_stage3}
    &q_{j}^{(0,3)} = \frac{1}{1 + \gamma\Phi}\left[1 - \tilde{a}_{31}\Phi - a_{31}\nu\frac{\exp\rpth{ik\Delta x} - \exp\rpth{-ik\Delta x}}{2} - \right. \nonumber \\
    &\hspace{2.6cm}\left. \tilde{a}_{32}\frac{\Phi}{2}G^{(n,2)}\rpth{\exp\rpth{-ik\Delta x} + 1} + a_{32}\nu G^{(n,2)}\rpth{\exp\rpth{-ik\Delta x} - 1}\right]q_{j}^{0}.
\end{align}
Denote by $G^{(n,3)}$
\begin{align}
    &G^{(n,3)} = \frac{1}{1 + \gamma\Phi}\left[1 - \tilde{a}_{31}\Phi - a_{31}\nu\frac{\exp\rpth{ik\Delta x} - \exp\rpth{-ik\Delta x}}{2} - \right. \nonumber \\
    &\hspace{2.7cm}\left. \tilde{a}_{32}\frac{\Phi}{2}G^{(n,2)}\rpth{\exp\rpth{-ik\Delta x} + 1} + a_{32}\nu G^{(n,2)}\rpth{\exp\rpth{-ik\Delta x} - 1}\right]
\end{align}
and notice that
\begin{equation}\label{eq:deltaq_von_Neumann_stage3}
    q_{j+1}^{(0,3)} - q_{j-1}^{(0,3)} = G^{(n,3)}\rpth{\exp\rpth{ik\Delta x} - \exp\rpth{-ik\Delta x}}q_{j}^{0}.
\end{equation}
Finally, considering the final update, one obtains
\begin{align}\label{eq:amplification_von_Neumann_tmp}
    &q_{j}^{1} = \left[\rpth{1 - \tilde{a}_{31}\Phi - \tilde{a}_{31}\nu\frac{\exp\rpth{ik\Delta x} - \exp\rpth{-ik\Delta x}}{2}}q_{j}^{0} - \right. \nonumber \\
    &\hspace{1.0cm}\left.\tilde{a}_{32}\rpth{\Phi\frac{q_{j-\frac{1}{2}}^{(0,2)} + q_{j+\frac{1}{2}}^{(0,2)}}{2} - \nu\rpth{q_{j-\frac{1}{2}}^{(0,2)} - q_{j+\frac{1}{2}}^{(0,2)}}} - \gamma\rpth{\Phi q_{j}^{(0,3)} + \nu\frac{q_{j+1}^{(0,3)} - q_{j-1}^{(0,3)}}{2}}\right]
\end{align}
Hence, substituting \eqref{eq:deltaq_von_Neumann_stage2}-\eqref{eq:sumq_von_Neumann_stage2} and \eqref{eq:deltaq_von_Neumann_stage3}, one obtains
\begin{align}
    &q_{j}^{1} = \left[\rpth{1 - \tilde{a}_{31}\Phi - \tilde{a}_{31}\nu\frac{\exp\rpth{ik\Delta x} - \exp\rpth{-ik\Delta x}}{2}} - \right. \nonumber \\
    &\hspace{1.0cm}\left.\tilde{a}_{32}\rpth{\Phi\frac{G^{(n,2)}\rpth{\exp\rpth{-ik\Delta x} + 1}}{2} - \nu\rpth{G^{(n,2)}\rpth{\exp\rpth{-ik\Delta x} + 1}}} - \right. \nonumber \\
    &\hspace{1.0cm}\left.\gamma\rpth{\Phi G^{(n,3)} + \nu\frac{G^{(n,3)}\rpth{\exp\rpth{ik\Delta x} - \exp\rpth{-ik\Delta x}}}{2}}\right]q_{j}^{0},
\end{align}
so that the amplification factor $G$ reads as follows
\begin{align}\label{eq:amplification_von_Neumann}
    &G = \left[\rpth{1 - \tilde{a}_{31}\Phi - \tilde{a}_{31}\nu\frac{\exp\rpth{ik\Delta x} - \exp\rpth{-ik\Delta x}}{2}} - \right. \nonumber \\
    &\hspace{1.0cm}\left.\tilde{a}_{32}\rpth{\Phi\frac{G^{(n,2)}\rpth{\exp\rpth{-ik\Delta x} + 1}}{2} - \nu\rpth{G^{(n,2)}\rpth{\exp\rpth{-ik\Delta x} + 1}}} - \right. \nonumber \\
    &\hspace{1.0cm}\left.\gamma\rpth{\Phi G^{(n,3)} + \nu\frac{G^{(n,3)}\rpth{\exp\rpth{ik\Delta x} - \exp\rpth{-ik\Delta x}}}{2}}\right].
\end{align}
One can verify that
\begin{equation}
    \lim_{\Phi \to \infty} G = G_{lim},
\end{equation}
as defined in \eqref{eq:Glim}. The reader can refer to the Mathematica~\cite{Mathematica2024} notebook \texttt{von\_Neumann.nb} in the GitHub repository \url{https://github.com/federicg/lava-flow.git} for the actual computation.

%%%%%%%%%% Bibliography %%%%%%%%%%%%%%%%
\printbibliography[category=cited]
\printbibliography[title={Further Reading}, notcategory=cited, resetnumbers=true]

\end{document}